\documentclass[12pt, eqno]{article}
\usepackage{bbm}
\usepackage{mathrsfs}
\usepackage{amsfonts}
\usepackage{amssymb}
\usepackage{graphicx}
\usepackage[all]{xy}

\usepackage{psfrag}
\usepackage{subfigure}
\usepackage{color}

\usepackage{amssymb,latexsym}
\usepackage{amsmath,latexsym}
\usepackage{amscd}

\newcommand{\N}{{\mathbb N}}

\newcommand{\R}{{\mathbb R}}

\newtheorem{theorem}{Theorem}[section]
\newtheorem{cor}[theorem]{Corollary}

\newtheorem{thm}[theorem]{Theorem}

\newtheorem{remark}[theorem]{Remark}
\newtheorem{rmk}[theorem]{Remark}
\newtheorem{lemma}[theorem]{Lemma}

\newtheorem{prop}[theorem]{Proposition}
\newtheorem{proposition}[theorem]{Proposition}
\newtheorem{claim}[theorem]{Claim}
\newtheorem{conjecture}{Conjecture}

\begin{document}

\title{A class of new bi-invariant metrics on\\ the Hamiltonian diffeomorphism groups}

\author{Guangcun Lu
\thanks{\endgraf
Partially supported
by the NNSF   10971014 and 11271044 of China,  PCSIRT, RFDPHEC (No. 200800270003)
 and the Fundamental Research Funds for the Central Universities (No. 2012CXQT09).}\hspace{10mm}Tie Sun\\
{\normalsize School of Mathematical Sciences, Beijing Normal University}\\
{\normalsize Laboratory of Mathematics and Complex Systems, Ministry of Education}\\
{\normalsize Beijing 100875, The People's Republic of China}\\
{\normalsize (gclu@bnu.edu.cn\hspace{5mm}suntie@mail.bnu.edu.cn)}}
\maketitle \vspace{-0.1in}

\abstract{In this paper, we construct infinitely many bi-invariant
metrics on the Hamiltonian diffeomorphism group and study their basic properties
and corresponding generalizations of the Hofer inequality and Sikorav one. }\\




\section{Introduction and Main Results}\setcounter{equation}{0}

\subsection{The Hofer metric}

In 1989, H. Hofer \cite{Ho90} constructed a remarkable  bi-invariant Finsler metric on
the group of compactly supported Hamiltonian diffeomorphisms ${\rm Ham}(M,\omega)$ of
a symplectic manifold $(M,\omega)$, nowadays known as Hofer metric.
Since then the intrinsic geometry of it has been being a very active and fruitful research
field in  symplectic topology and Hamiltonian dynamics (see the books
\cite{HoZe, McSa, PoBook}, and the surveys \cite{En14, Mc10, Po98ICM, Os14}
 and references therein for current progress situation).

Especially, a recent celebrated result  made by Buhovsky and  Ostrover \cite{BuOs11}
is a positive answer to the uniqueness question of the Hofer metric raised by
Eliashberg and Polterovich \cite{EP93}. They showed that up to equivalence of metrics
the Hofer metric is the only bi-invariant Finsler metric on the group of Hamiltonian diffeomorphisms
of a closed symplectic manifold under a natural assumption. For studies of non-Finslerian
bi-invariant  metrics on ${\rm Ham}(M,\omega)$ the readers may refer to \cite{Vi92, Sc, Oh05}.

Let us briefly review the construction of the Hofer metric
following the notations in \cite{PoBook} without special statements.
The readers who are familiar with it may directly read the next section.
 Let $(M^{2n},\omega)$ be a connected
symplectic manifold of dimension $2n$ without boundary.  Denote by $\mathcal{A}(M)$
the space of all smooth functions on $M$ with compact support (resp.  zero-mean with
respect to the canonical volume form $\omega^n$) if $M$ is open (resp. closed).
A (time-dependent) smooth Hamiltonian function $F$ on $M\times I$,
where $I\subset\R$ is an interval, is called \textbf{normalized} if
$F_t=F(\cdot,t)$ belongs to $\mathcal{A}(M)$ for all $t$, and
$\cup_{t\in I}{\rm supp}(F_t)$ is contained in a compact subset of
$M$ in the case when $M$ is open. Such a normalized $F$ determines a
(time-dependent) Hamiltonian vector field $X_{F_t}$ on $M$ via
$i_{X_{F_t}}\omega=-dF_t$,  and when $I=[0,1]$ the corresponding
flow $\{f_t\}$ starting from the identity is called a
\textbf{Hamiltonian isotopy} generated by $F$ and is also denoted by
$\{\phi_F^t\}$ for convenience.
 A diffeomorphism of $M$ is said to be \textbf{Hamiltonian} if it
can be represented as a time-one map of some Hamiltonian isotopy.
Denote by ${\rm Ham}(M, \omega)$ the set of all Hamiltonian
diffeomorphisms on $(M, \omega)$. It is  a subgroup of ${\rm
Symp}(M, \omega)$. The Lie algebra of ${\rm Ham}(M, \omega)$ can be
naturally identified with $\mathcal{A}(M)$ because of the following important fact by Banyaga
(cf. \cite[Prop.1.4.B]{PoBook}).

\begin{lemma}[\hbox{\cite{Ban78}}]\label{lem:Ban}
For every smooth path $\{f_t\}$ in ${\rm Ham}(M, \omega)$,
$t\in[a,b]$, there exists a unique (time-dependent) normalized Hamiltonian
function $F : M\times[a,b]\rightarrow \mathbb{R}$ such that
\begin{equation}\label{e:ode}
\frac{d}{dt}f_t(x)=X_{F_t}(f_t(x))\quad\forall (x,t)\in M\times [a,b].
\end{equation}
\end{lemma}
Hence the tangent
vector of the path $\{f_t\}$ at $t=s$ is the function $F_s$. The
adjoint action of ${\rm Ham}(M, \omega)$ on $\mathcal{A}(M)$ is
given by ${\rm Ad}_fG=G\circ f^{-1}$ for $f\in {\rm Ham}(M, \omega)$
and $G\in\mathcal{A}(M)$. Any adjoint invariant norm
$\|\cdot\|$ on $\mathcal{A}(M)$, i.e., $\|{\rm Ad}_fG\|=\|G\circ
f^{-1}\|=\|G\|$ for any $f\in {\rm Ham}(M, \omega)$ and
$G\in\mathcal{A}(M)$, defines a Finsler structure on ${\rm Ham}(M,
\omega)$, and thus  the length of a Hamiltonian path
$\{f_t\}$, $t\in [a,b]$ with (unique) normalized Hamiltonian $F$ by
\begin{equation}\label{def:length}
{\rm Length}\{f_t\}=\int_a^b\|F_t\|dt,
\end{equation}
which does not depend on the parametrization.
Without loss of generality we could fix $a=0,b=1$ in the definition.
For arbitrary $\phi,\varphi\in {\rm Ham}(M, \omega)$, their
pseudo-distance is defined by
\begin{equation}\label{e:pseudo-dis}
d(\phi,\varphi)=\inf\{{\rm Length}(\alpha)\}
\end{equation}
where the infimum is taken over all smooth Hamiltonian path
$\alpha:[a,b]\rightarrow\mathcal {\rm Ham}(M, \omega)$ with
$\alpha(a)=\phi$ and $\alpha(b)=\varphi$. It is a bi-invariant
pseudo-metric.  (If the norm $\|\cdot\|$ on $\mathcal{A}(M)$ is not adjoint
invariant, $d$ is only right-invariant.)   The non-degeneracy axiom,
$d(\phi,\varphi)>0$ for $\phi\neq\varphi$, is not satisfied in
general.

When the norm $\|\cdot\|$ is chosen as the $L_\infty$-norm,
$$
\|H\|_\infty :=\max_xH-\min_xH,\quad \hbox{for}\quad
H\in\mathcal{A}(M),
$$
Hofer  showed  in \cite{Ho90} that the corresponding pseudo-metric
$d_H$ is  a genuine metric in the case $M=\mathbb{R}^{2n}$. Later,
this result was generalized   to some larger class of symplectic
manifolds by Polterovich  \cite{Po93}, and finally to general
manifolds by Lalonde and McDuff  \cite{LaMc}. Nowadays  this
bi-invariant Finsler metric on the group ${\rm Ham}(M, \omega)$ is called the
{\bf Hofer metric}, and the function $\|\cdot\|_H=d_H(\cdot,{\rm id}_M):{\rm Ham}(M, \omega)\to\mathbb{R}$ is called the {\bf Hofer norm}.
Let
\begin{eqnarray*}
&&{\mathcal H}={\mathcal H}(M)=\{H\in C^\infty(M\times [0,1])\,|\, H_t=H(\cdot,t)\in{\mathcal A}(M)\;\forall t\in [0,1]\},\\
&&{\mathcal F}={\mathcal F}(M)=\{H\in C^\infty(M\times\mathbb{R}/\mathbb{Z})\,|\, H_t=H(\cdot,t)\in{\mathcal A}(M)\;\forall t\in \mathbb{R}\}.
\end{eqnarray*}
Every $\phi\in{\rm Ham}(M,\omega)$ can be written as $\phi_F^1$ with $F\in\mathcal{F}$. Moreover it holds that
\begin{eqnarray}\label{e:carse}
&&d_H(\phi,\varphi)=\inf\big\{\int^1_0\|H_t-K_t\|_\infty dt\,\big|\,H\in{\mathcal H}\;\hbox{generates}\;
\varphi\;\&\;K\in{\mathcal H}\;
\hbox{generates}\;\phi\big\},\nonumber\\
&&d_H({\rm id}_M,\varphi)=\inf\big\{\max_t\|F_t\|_\infty \,\big|\,F\in{\mathcal F}\;\hbox{generates}\;
\varphi\big\},
\end{eqnarray}
see \cite[(5.9)]{HoZe} for the first one, and \cite[Lemma 5.1.C]{PoBook} for the second.

\subsection{New bi-invariant metrics}

Our new  bi-invariant metrics on ${\rm Ham}(M,\omega)$ will be constructed
in a similar way to Hofer's. For a
smooth path $f:[0,1]\to {\rm Ham}(M,\omega)$  generated by a
(time-dependent) normalized Hamiltonian function $F$, and each
integer $k=0,1,2\cdots$, we define the \textbf{k-length} of $f$
by
$$
{\rm Length}_k(f):=\sum_{i=0}^{k}\int_0^1\left\|
\frac{\partial^i F_t}{\partial t^i}\right\|_\infty dt=
\sum_{i=0}^{k}\int_0^1\Big(\max_x
\frac{\partial^i F_t}{\partial t^i}-\min_x\frac{\partial^i F_t}
{\partial t^i}\Big)dt.
$$
Clearly, ${\rm Length}_0$ is the same as that of (\ref{def:length}).
However, unlike ${\rm Length}_0$ the $k$-length ($k\ge 1$)
strongly depends on the choice of parametrization.

Call  a continuous path $f:[0, 1]\to {\rm Ham}(M,\omega)$  {\bf piecewise smooth} if
there exists a division $0=t_0<t_1<\cdots<t_n=1,n\in\mathbb{N}$, such that
$f_i=f|_{[t_{i-1}, t_i]}$ is smooth for $i=1,\cdots,n$.
Let  $F=\{F^l\}_{l=1}^n$ be the corresponding normalized Hamiltonian function,
that is,
$$\frac{d}{dt}f_t(x)=X_{F^i_t}(f_t(x))\quad\forall (x,t)\in M\times [t_{i-1}, t_i],
i=1,\cdots,n.
$$
Define its $k$-length by
\begin{equation}\label{e:length_k}
{\rm Length}_k(f):=\sum^n_{i=1}{\rm Length}_k(f_i)=\sum^n_{l=1}\sum_{i=0}^{k}\int_{t_{l-1}}^{t_l}\left\|
\frac{\partial^i F^l_t}{\partial t^i}\right\|_\infty dt.
\end{equation}
For $\phi,\varphi\in {\rm Ham}(M,\omega)$,
let $\Omega(\phi,\varphi)$
denote the space of all continuous and piecewise smooth paths
$f:[0,1]\rightarrow {\rm Ham}(M,\omega)$ from $f(0)=\phi$ to
$f(1)=\varphi$.  Then we define
pseudo-distances between $\phi$ and $\varphi$ by
\begin{equation}\label{def:k-dis}
d_k(\phi,\varphi)=\inf\left\{{\rm
Length}_k(f)\,\Big|\,f\in\Omega(\phi,\varphi)\right\},\quad
k=0,1,\cdots.
\end{equation}

\begin{lemma}\label{lem:HoferM}
$d_0=d_H$ on ${\rm Ham}(M,\omega)$.
\end{lemma}

 Clearly, $d_0\le d_H$ since (\ref{def:length})
 does not depend on the parametrization. The converse inequality easily follows
 from the definition of $d_0$ and the triangle inequality for $d_H$.

Let us make some comments on the definition of $d_k$.

\begin{rmk}
{\rm (i) For a positive integer $k$,  the function
${\rm Length}_k(\cdot)$ depends on the choice of parametrization of the path,
and the derivative terms will in fact vanish when we take the infimum with respect
to the path space with variant parametrization intervals. In fact, suppose that a smooth Hamiltonian path
$\alpha:[0,1]\to {\rm Ham}(M,\omega)$ is generated by a normalized Hamiltonian
$F_t$.
For any $b>0$,  the reparametrized path
$$
\beta_b(t):=\alpha(t/b):[0,b]\to {\rm Ham}(M,\omega)
$$
is generated by the Hamiltonian function $G(x,t)=\frac{1}{b}
F(x,\frac{t}{b})$, and hence
$$
\int_0^b\Big\|\frac{\partial G_t}{\partial t}\Big\|_\infty dt=
\int_0^b\frac{1}{b^2}\Big\|\frac{\partial F}{\partial t}(x,\frac{s}{b})\Big\|
_\infty ds=
\frac{1}{b}\int_0^1\Big\|\frac{\partial F}{\partial t}(x,t)\Big\|_\infty dt\to 0
$$
as $b\rightarrow+\infty$.  This fact still holds for the higher order derivatives. It follows that
$$
b\to +\infty \Longrightarrow {\rm Length}_k(\beta_b) \rightarrow{\rm
Length}_0(\alpha).
$$
Thus if we define $d_k$ as in (\ref{e:pseudo-dis}) no new thing can be obtained. This is why
we fix the parametrization interval of paths, $[a,b]=[0,1]$. \\
(ii) Since the definition of $k$-length involves  the derivatives of a
Hamiltonian function until $k$ order the quasi-triangle inequality cannot be obtained if
we restrict to smooth paths from $[0,1]$ to ${\rm Ham}(M,\omega)$ in the definition of
$d_k$. It is this reason that  we extend the space of all
smooth Hamiltonian paths to include the piecewise smooth ones.}
\end{rmk}

As expected $d_k$  has the following properties.

\begin{theorem}\label{th:property}
$d_H=d_0\le d_1\le d_2\le\cdots$, and
\begin{description}
\item[(i)] {\rm (symmetry)} $d_k(\phi,\varphi)=d_k(\varphi,\phi)$,
\item[(ii)]{\rm (quasi-triangle inequality)} $d_k(\phi,\varphi)\leq 2^k
    (d_k(\phi,\theta)+d_k(\theta,\varphi))$,
\item[(iii)] {\rm (non-degeneracy)} $d_k(\phi,\varphi)\geq0$, and
    $d_k(\phi,\varphi)=0\Longleftrightarrow\phi=\varphi$,
\item[(iv)] {\rm (bi-invariance)} $d_k(\phi,\varphi)=d_k(\phi\theta,\varphi\theta)
    =d_k(\theta\phi,\theta\varphi)$,
\end{description}
for any $\phi,\varphi,\theta\in{\rm Ham}(M,\omega)$ and
$k=0,1,\cdots$.
\end{theorem}

This theorem shows that  (\ref{def:k-dis}) gives  a
sequence of bi-invariant  quasi-metrics $\{d_k\}^\infty_{k=0}$ on ${\rm Ham}(M,\omega)$.
Recall that  a {\bf quasidistance} on a nonempty set $X$ is a function
 $\rho : X\times X\rightarrow[0,+\infty)$
such that (i) $\rho(x,y)=\rho(y,x)$ for all $x,y\in X$,
(ii) $\rho(x,y)=0$ if and only if $x=y$,
(iii)  and there exists a finite constant
$c\geq1$ (quasi-triangle constant) such that
$\rho(x,y)\leq c(\rho(x,z)+\rho(z,y))$  for every $x,y,z\in X$.
  Such a pair $(X, \rho)$ is called a {\bf quasimetric} space.
    (See \cite{MiMMM}). A {\bf group norm} (resp. {\bf quasinorm}) on a group $G$
is a symmetric, nondegenerate and nonnegative function $\psi$ which
is subadditive (resp. $C$-subadditive  for some finite constant $C\ge 1$, that is,
$\psi(xy)\le C(\psi(x)+\psi(y))$ for all $x,y\in G$). (See \cite[page 113]{MiMMM}).

  For every $d_k$, let us define a function
\begin{equation}\label{e:quasi-norm}
\|\cdot\|_k=d_k(\cdot,{\rm id}_M).
\end{equation}
From Theorem~\ref{th:property} one easily derives:

\begin{theorem}\label{th:quasi-norm}
$\|\cdot\|_0=\|\cdot\|_H$, and for every $k\in\N$, $\|\cdot\|_k$ is a quasinorm,
precisely speaking it satisfies:
\begin{description}
\item[(i)] {\rm (Symmetry)} $\|\phi^{-1}\|_k=\|\phi\|_k$,
\item[(ii)]{\rm (The quasi-triangle inequality)} $\|\phi\varphi\|_k\leq
    2^k(\|\phi\|_k+ \|\varphi\|_k)$,
\item[(iii)] {\rm (Non-degeneracy)} $\|\phi\|_k\geq0$, and $\|\phi\|_k
    =0\Longleftrightarrow\phi={\rm id}_M$,
\item[(iv)] {\rm (Conjugate invariance)}
    $\|\theta\varphi\theta^{-1}\|_k=\|\varphi\|_k$,
\end{description}
where $\phi,\varphi,\theta\in{\rm Ham}(M,\omega)$ are arbitrary.
\end{theorem}

This  shows that every $\|\cdot\|_k$ is a conjugate invariant
quasinorm on the group ${\rm Ham}(M,\omega)$.
Apply Theorem~\ref{th:MiMMM} to $G={\rm Ham}(M,\omega)$ and
$\psi=\|\cdot\|_k$, $k\in\N$, we get

\begin{thm}\label{th:MiMMM-d_k}
Define the function $\|\!|\cdot\|\!|_k:{\rm Ham}(M,\omega)\to [0, \infty)$ by
\begin{equation}\label{e:MiMMM1-d_k}
\|\!|\phi\|\!|_k:=\inf\left\{\left(\sum^N_{i=1}\|\phi_i\|_k^{\frac{1}{k+1}}\right)^
{1+k}\Bigg|\;\begin{array}{ll} &N\in\N,\; (\phi_1,\cdots,\phi_N)\in
{\rm Ham}(M,\omega)^{(N)},\\
&\phi=\phi_1\cdots \phi_N
\end{array}
\right\},
\end{equation}
where ${\rm Ham}(M,\omega)^{(N)}=\underbrace{{\rm
Ham}(M,\omega)\times\cdots\times {\rm
Ham}(M,\omega)}_{N\hbox{¸ö}}$, $N\in\mathbb{N}$. Then
\begin{equation}\label{e:MiMMM2-d_k}
\|\!|\phi\|\!|_k=\|\!|\phi^{-1}\|\!|_k,\quad\forall \phi\in{\rm
Ham}(M,\omega),
\end{equation}
\begin{equation}\label{e:MiMMM3-d_k}
\|\!|\theta\phi\theta^{-1}\|\!|_k=\|\!|\phi\|\!|_k,\quad\forall
\theta, \phi\in{\rm Ham}(M,\omega),
\end{equation}
\begin{equation}\label{e:MiMMM4-d_k}
4^{-(1+k)}\|\phi\|_k\le\|\!|\phi\|\!|_k\le\|\phi\|_k,\quad\forall
\phi\in{\rm Ham}(M,\omega),
\end{equation}
and
\begin{equation}\label{e:MiMMM5-d_k}
\|\!|\phi\varphi\|\!|_k^\beta\le\|\!|\phi\|\!|_k^\beta+
\|\!|\varphi\|\!|_k^\beta,\quad\forall \phi, \varphi\in{\rm
Ham}(M,\omega)
\end{equation}
for each $\beta\in (0, \frac{1}{1+k}]$.
Also, for each $N\in\N$, $\beta\in (0,\frac{1}{1+k}]$,
and $\phi_i\in{\rm Ham}(M,\omega)$, $i=1,\cdots,N$,  it holds that
\begin{equation}\label{e:MiMMM6-d_k}
\|\phi_1\cdots \phi_N\|_k \le
4^{k+1}\left\{\sum^N_{i=1}\|\phi_i\|_k^\beta\right\}^{\frac{1}{\beta}},
\end{equation}
and hence for each sequence $(\phi_i)_{i\in\N}\subset {\rm Ham}(M,\omega)$,
\begin{equation}\label{e:MiMMM7-d_k}
\sup_{N\in\N}\|\phi_1\cdots\phi_N\|_k \le
4^{k+1}\left\{\sum^\infty_{i=1}\|\phi_i\|_k^\beta\right\}^{\frac{1}{\beta}}.
\end{equation}
\end{thm}

(\ref{e:MiMMM2-d_k}) can be derived from (\ref{e:MiMMM1-d_k}) and
the symmetry of $\|\cdot\|_k$, and
(\ref{e:MiMMM3-d_k}) can be obtained by Remark~\ref{rm:MiMMM} and the
conjugate invariance of $\|\cdot\|_k$.

\begin{cor}\label{cor:MiMMM-d_k}
The function $\|\!|\cdot\|\!|_k$ defined in (\ref{e:MiMMM1-d_k}) is
a conjugate invariant quasi-norm on ${\rm Ham}(M,\omega)$ which is
equivalent to $\|\cdot\|_k$; and for every $\beta\in (0, \frac{1}{1+k}]$,
$\|\!|\cdot\|\!|_k^\beta$ is a conjugate invariant norm on ${\rm
Ham}(M,\omega)$. Thus
\begin{equation}\label{e:MiMMM8-d_k}
{\tilde d}_{k}(\phi,\varphi):=\|\!|\phi\varphi^{-1}\|\!|_k
\end{equation}
is a bi-invariant quasimetric on ${\rm Ham}(M,\omega)$; and for each $\beta\in (0,
\frac{1}{1+k}]$,
\begin{equation}\label{e:MiMMM9-d_k}
\tilde{d}^\beta_{k}(\phi,\varphi):=
({\tilde d}_{k}(\phi,\varphi))^\beta
\end{equation}
is a bi-invariant metric on ${\rm Ham}(M,\omega)$.
They all induce the same topology as $d_k$.
\end{cor}

Consider the commutator  of
two elements $\varphi$ and $\psi$ in ${\rm Ham}(M,\omega)$, $[\varphi,\psi]:=\varphi\psi\varphi^{-1}\psi^{-1}$.
It follows from Theorem~\ref{th:quasi-norm} that
\begin{equation}\label{e:commut1}
\|[\varphi,\psi]\|_k\le 2^{k+1}\min\{\|\varphi\|_k,\|\psi\|_k\}.
\end{equation}
Similarly (\ref{e:MiMMM2-d_k})-(\ref{e:MiMMM5-d_k}) lead to
\begin{equation}\label{e:commut2}
\|\!|[\varphi,\psi]\|\!|_k^\beta\le 2\min\{\|\!|\varphi\|\!|_k^\beta, \|\!|\psi\|\!|_k^\beta\}
\end{equation}
for each $\beta\in (0, \frac{1}{1+k}]$. For a non-empty subset $A\subset M$ let
$e_k(A)$ (resp. $\tilde{e}_k(A)$) denote the {\bf displacement energy} of it with respect to $\|\cdot\|_k$
(resp. $\|\!|\cdot\|\!|_k$), that is,
\begin{eqnarray}\label{e:disE1}
&&e_k(A)=\inf\{\|\vartheta\|_k\,\big|\, \vartheta\in{\rm Ham}(M,\omega)\;\&\;A\cap\vartheta(A)=\emptyset\},\\
&&\tilde{e}_k(A)=\inf\{\|\!|\vartheta\|\!|_k\,\big|\, \vartheta\in{\rm Ham}(M,\omega)\;\&\;A\cap\vartheta(A)=\emptyset\}.\label{e:disE2}
\end{eqnarray}
As in the proof of \cite[Lemma~2.3.B]{EP93} we may obtain:

\begin{theorem}\label{th:e-estimate}
Let $U\subset M$ be a non-empty open subset. Then for any $\varphi,\psi\in {\rm Ham}(M,\omega)$
with ${\rm supp}(\varphi)\subset U$ and ${\rm supp}(\psi)\subset U$ it holds that
\begin{eqnarray}\label{e:e-e}
\|[\varphi,\psi]\|_k\le 4^{k+1}e_k(U)\quad\hbox{and}\quad
\|\!|[\varphi,\psi]\|\!|_k^\beta\le 4(\tilde{e}_k(U))^\beta
\end{eqnarray}
for each $\beta\in (0, \frac{1}{1+k}]$.
\end{theorem}

Motivated by the so-called ``coarse" Hofer norm, for  $f=\{f_t\}\in\Omega(\phi,\psi)$,
if $F=\{F^l\}_{l=1}^n$ is the corresponding normalized Hamiltonian function,
we use
\begin{equation}\label{def:k-length*}
{\rm Length}^\ast_k(f):=|||F|||_k:=\sum_{i=0}^k\max_{1\le l\le n}\max_{t_{l-1}\le t\le t_l}\left\|\frac{\partial^iF^l}{\partial t^i}(x,t)\right\|_\infty
\end{equation}
(which is independent of the choices of divisions) to  replace (\ref{e:length_k}), and  obtain another sequence of bi-invariant quasimetrics
\begin{equation}\label{def:k-dis*}
d^\ast_k(\phi,\varphi)=\inf\left\{{\rm
Length}^\ast_k(f)\,\Big|\, f\in\Omega(\phi,\varphi)\right\},\quad
k=0,1,\cdots
\end{equation}
as in (\ref{def:k-dis}).  Clearly, $d_k\le d^\ast_k$ for any $k\in\mathbb{N}\cup\{0\}$.
(\ref{e:carse}) implies  $d^\ast_0=d_0=d_H$. Let $\|\cdot\|^\ast_k=d_k^\ast(\cdot,{\rm id}_M)$.
Correspondingly, we have also $\|\!|\cdot\|\!|^\ast_k$ as in (\ref{e:MiMMM1-d_k}), and
$\tilde{d}^\ast_{k}$, $e_k^\ast$ and $\tilde{e}_k^\ast$.


\begin{theorem}\label{th:corres}
All conclusions from Theorem~\ref{th:property} to Theorem~\ref{th:e-estimate}
still hold for $\|\cdot\|^\ast_k$, $d_k^\ast$, $\|\!|\cdot\|\!|^\ast_k$ and
$\tilde{d}^\ast_{k}$, $e_k^\ast$ and $\tilde{e}_k^\ast$, but we need to add a factor $2$ for
the coefficients of inequalities in (ii) of Theorems~\ref{th:property},~\ref{th:quasi-norm}, (\ref{e:commut1})
and the factor $4$ in the first inequality of (\ref{e:e-e}).
\end{theorem}

The following result shows that the Hofer inequality in \cite{Ho93} also holds for each $d_k^\ast$.

\begin{theorem}\label{th:HoferE}
For every $\varphi,\psi\in{\rm Ham}(\mathbb{R}^{2n},\omega_0)$,
$$
d_k^\ast(\varphi,\psi):=\|\varphi\psi^{-1}\|_k^\ast\leq C\,{\rm diameter}\bigl({\rm supp}(\varphi\psi^{-1})\bigr)
|\varphi-\psi|_{C^0},
$$
where $C$ is a constant and $C\leq2^{3k+8}(k+1)^2\big(1+2^{k+1}+2^{2k+2}+2^{3k+3}\big)$.
({\rm Note}: if $k=0$ the constant $C$ can be chosen as $128$ as in the Hofer inequality.)
\end{theorem}

Similarly,  for any subset $S\subset\mathbb{R}^{2n}$, we define the {\bf coarse proper displacement
$k$-energy} $e_{p,k}^\ast(S)$ of it as
\begin{eqnarray*}
e_{p,k}^\ast(S)=\bigl\{a>0&\big|&\hbox{for every bounded subset}\,A\subset S\;
\exists\,\psi\in{\rm Ham}(\mathbb{R}^{2n},\omega_0)\\
&&\hbox{such that}\,
\|\psi\|_k^\ast\leq a\,\hbox{and}\,A\,\hbox{and}\,\psi(A)\,\hbox{are properly separated}\bigr\},
\end{eqnarray*}
and get the following generalization of the Sikorav inequality.

\begin{theorem}\label{th:SikE}
If $H\in\mathcal{H}(\mathbb{R}^{2n})$ satisfies ${\rm supp}(H)\subset U\times[0,1]$, then
$$
\|\varphi_H\|_k^\ast\leq2^{2k+4}(k+1)\big(1+2^{k+1}+2^{2k+2}+2^{3k+3}\big)e_{p,k}^\ast(U).
$$
\end{theorem}

Finally, let us discuss the corresponding question investigated by Eliashberg and Polterovich \cite{EP93}.
When $0<p<\infty$ the function
$$
\mathcal{A}(M)\ni H\to \|H\|_p=\left(\int_M|H|^p\omega^n\right)^{\frac{1}{p}}
$$
is an adjoint invariant quasinorm because
\begin{equation}\label{e:adjointAction}
\int_M|G(f^{-1}(x))|^p\omega^n=\int_M(f^{-1})
^*(|G(x)|^p\omega^n)=\int_M|G(x)|^p\omega^n
\end{equation}
for any $G\in \mathcal{A}(M)$ and $f\in{\rm Symp}(M,\omega)$, and
$$\|H+G\|_p\le
K_p(\|H\|_p+\|G\|_p),\quad\forall H, G\in
\mathcal{A}(M),
$$
where  $K_p$ is equal to $1$ for $p\ge 1$, and $2^{\frac{1-p}{p}}$ for $0<p<1$ (so
 $\|\cdot\|_p$ is only a quasinorm on $\mathcal{A}(M)$ in this case).

For a smooth path $f:[0,1]\to {\rm Ham}(M,\omega)$  generated
by a (time-dependent) normalized Hamiltonian function $F_t$, and
each integer $k=0,1,2\cdots$, we define the \textbf{(k,p)-length} of
$f$ by
\begin{equation}\label{e:length_k-p}
{\rm Length}_{(k,p)}(f):=\sum_{i=0}^{k}\int_0^1\left\|
\frac{\partial^i F_t}{\partial t^i}\right\|_pdt,
\end{equation}
and the \textbf{(k,p)-length} of $f\in\Omega(\phi,\varphi)$ by
the sum of $(k,p)$-lengths of all smooth pieces of it.
By the same proofs as those of Theorem~\ref{th:property}
it is readily verified that
\begin{equation}\label{def:(k,p)-dis}
d_{(k,p)}(\phi,\varphi)=\inf\left\{{\rm
Length}_{(k,p)}(f)\,\Big|\,f\in\Omega(\phi,\varphi)\right\}
\end{equation}
defines a pseudo quasimetric $d_{(k,p)}$ on ${\rm Ham}(M,\omega)$ for any
$p>0, k=0,1,\cdots$. Note that
\begin{equation}\label{def:(k,p)-dis1}
d_{(0,p)}(\phi,\varphi)=\inf\left\{{\rm
Length}_{(0,p)}(f)\,\Big|\,f\in\Omega(\phi,\varphi)\;\hbox{is smooth}\right\}.
\end{equation}
Eliashberg and Polterovich \cite{EP93, PoBook} showed for each $p\in [1, \infty)$
that the pseudo-distance $d_{(0,p)}$ is degenerate, and
 vanishes if $M$ is
closed. We have the following extension.


\begin{theorem}\label{th:EP}
For each $k\in\N\cup\{0\}$ and each $0<p<1/k$, the pseudo quasimetric $d_{(k,p)}$ is degenerate,
and vanishes if $M$ is closed.
\end{theorem}

Similarly, for a closed embedded Lagrangian submanifold $L$ of $(M, \omega)$
 let $\mathcal{L}(M,\omega, L)$ denote the space of Lagrangian submanifolds of
$(M, \omega)$ which is Hamiltonian isotopic to $L$. For each $k=0,1,2,\cdots$,
define $\delta_k:\mathcal{L}(M,\omega, L)\times\mathcal{L}(M,\omega, L)\to\mathbb{R}\cup\{+\infty\}$
by
$$
\delta_k(L_1,L_2)=\inf\{\|\phi\|_k\,|\, \phi\in{\rm Ham}(M,\omega)\;\&\;\phi(L_1)=L_2\}.
$$
Then $\delta_k(L_1,L_2)=\delta_k(L_2,L_1)$ and $\delta_k(L_1,L_2)\le 2^k\bigl(\delta_k(L_1,L_3)+
\delta_k(L_2, L_3)\bigr)$ for any $L_i\in{\mathcal L}$, $i=1,2,3$. If $(M,\omega)$ is
a tame symplectic manifold, Chekanov showed in \cite{Che} that $\delta_0=\delta_H$ is non-degenerate,
and so each $\delta_k$ is a ${\rm Ham}(M,\omega)$-invariant quasimetric.

\

The paper is organized as follows. In Section 2 we give proofs of
Theorem~\ref{th:property},~\ref{th:HoferE}, ~\ref{th:SikE},
~\ref{th:EP}. Extensions of our metrics onto
the group of symplectic diffeomorphisms will be discussed in Section 3.
Finally, Section 4 outlines our constructions on
the group of strictly contact diffeomorphisms as a concluding remark.

\

\noindent{\bf Acknowledgements}. We would like to thank Professors
L. Polterovich, H.Hofer and J. C. Sikorav for  kind helps in understanding
their papers.

\section{Proofs}\setcounter{equation}{0}

\subsection{Proof of Theorem~\ref{th:property}}
In the following we always assume $k>0$.\\
(i) For $\Omega(\phi,\varphi)\ni f:[0,1]\rightarrow{\rm Ham}(M,\omega)$,
let $\bar{f}\in\Omega(\varphi,\phi)$ be defined by $\bar{f}(t)=f(1-t),\,\forall t\in
[0,1]$. Then $\Omega(\phi,\varphi)\ni f\to\bar{f}\in\Omega(\varphi,\phi)$
 is a bijection. By the assumption there exists
a division $0=t_0<t_1<\cdots<t_n=1, n\in\mathbb{N}$, such that
$f_i=f|_{[t_{i-1},t_i]}$ is smooth for $i=1,\cdots,n$. By
Lemma~\ref{lem:Ban}, for each $j=1,\cdots,n$
 there exists a (time-dependent) normalized Hamiltonian
function $F_j : M\times [t_{j-1},t_j]\rightarrow \mathbb{R}$ such
that
\begin{equation}\label{e:ode-j}
\frac{d}{dt}f_j(t)(x)=X_{F_j}\bigl(f_j(t)(x), t\bigr)
\end{equation}
for all $x\in M$ and $t\in [t_{j-1},t_j]$. Set $s_i=1-t_{n-i}$,
$i=0,\cdots,n$. Then $0=s_0<s_1<\cdots<s_n=1$ is a division of
$[0,1]$, and for any $t\in [s_{j-1}, s_j]$ it holds that
$$
\bar{f}_j(t)=\bar{f}|_{[s_{j-1}, s_j]}(t)=f(1-t)=f|_{[t_{n-j},
t_{n-j+1}]}(1-t)=f_{n-j+1}(1-t).
$$
From this and (\ref{e:ode-j}) it follows that
\begin{eqnarray*}
\frac{d}{dt}\bar{f}_j(t)(x)&=&-\frac{d}{ds}f_{n-j+1}(s)(x)|_{s=1-t}\\
&=&-X_{F_{n-j+1}}\bigl(f_{n-j+1}(s)(x),s\bigr)|_{s=1-t}=X_{G_j}\bigl(\bar{f}_j(t)(x),t\bigr),
\end{eqnarray*}
where $G_j:M\times [s_{j-1}, s_j]\to\R$ is given by
$G_j(x,t)=-F_{n-j+1}(x,1-t)$. Hence
$$
\frac{\partial^i (G_j(x, s))}{\partial
s^i}=(-1)^{i+1}\times\,(F_{n-j+1})^{(i)}_2(x,1-s)\quad\forall s\in
[s_{j-1}, s_j]
$$
for $i=0,\cdots,k$, where $(F_{n-j+1})^{(i)}_2$ stands for the $i$th
partial derivative of $F_{n-j+1}$ with respect to the second
variable. By the definition we have
\begin{eqnarray*}
\begin{split}
{\rm Length}_k(\bar{f})&=\sum^n_{j=1}{\rm Length}_k(\bar{f}_j)\\
&=\sum^n_{j=1}\sum_{i=0}^{k}\int_{s_{j-1}}^{s_j}\left\|
\frac{\partial^i (G_j(x, s))}{\partial s^i}\right\|_\infty ds \\
&=\sum^n_{j=1}\sum_{i=0}^{k}\int_{s_{j-1}}^{s_j}\left\|(-1)^{i+1}
(F_{n-j+1})^{(i)}_2(x,1-s)\right\|_\infty ds \\
&=\sum^n_{j=1}\sum_{i=0}^{k}\int_{s_{j-1}}^{s_j}\left\|(F_{n-j+1})
^{(i)}_2(x,1-s)\right\|_\infty ds \\
&=\sum^n_{j=1}\sum_{i=0}^{k}\int_{t_{n-j}}^{t_{n-j+1}}\left\|(F_{n-j+1})
^{(i)}_2(x,t)\right\|_\infty dt\\
&=\sum^n_{j=1}\sum_{i=0}^{k}\int_{t_{j-1}}^{t_{j}}\left\|(F_j)^{(i)}_2(x,t)\right\|_\infty
dt\\
&=\sum^n_{j=1}{\rm Length}_k({f}_j)={\rm Length}_k({f}).
\end{split}
\end{eqnarray*}
 Thus $d_k(\phi,\varphi)=d_k(\varphi,\phi)$.

\noindent{(ii)} Let $\Omega(\phi,\theta)\ni f:[0,1]\rightarrow{\rm Ham}(M,\omega)$,
$\Omega(\theta,\varphi)\ni g:[0,1]\rightarrow{\rm Ham}(M,\omega)$. We define
the product path $g\sharp f:[0,1]\rightarrow{\rm Ham}(M,\omega)$ of $f$ and $g$ by
\begin{eqnarray*}
g\sharp f(t)=
\begin{cases}
f(2t) & 0\leq t\leq 1/2\\
g(2t-1) & 1/2\leq t\leq 1,
\end{cases}
\end{eqnarray*}
then $g\sharp f$ is a piecewise smooth Hamiltonian path connecting $\phi$ and
$\varphi$, i.e. $g\sharp f\in\Omega(\phi,\varphi)$.

By assumption there exist divisions
\begin{center}
$0<t_0<t_1<\cdots<t_n=1$ and\;$0<t_0'<t_1'<\cdots<t_m'=1$,
\end{center}
such that for $i=1,\cdots,n$, $f_i=f|_{[t_{i-1},t_i]}$ is smooth; for
$j=1,\cdots,m$, $g_j=g|_{[t_{j-1}',t_j']}$ is also smooth. Denote the
Hamiltonian functions generating $\{f_i\}_{i=1}^n$, $\{g_j\}_{j=1}^m$ by
$\{F_i\}_{i=1}^n$ and $\{G_j\}_{j=1}^m$ respectively. For $0\leq i\leq n$,
set $s_i=\frac{t_i}{2}$, for $n+1\leq i\leq n+m$, set $s_i=\frac{t_{i-n}'+1}{2}$, then
$$
0=s_0<s_1<\cdots<s_n=\frac{1}{2}<s_{n+1}<\cdots<s_{n+m}=1
$$
is a division of $[0,1]$, and at this time $(g\sharp f)_i=g\sharp f|_{[s_{i-1},s_i]}$
is smooth for $i=1,\cdots,n+m$. Denote the Hamiltonian function generating
$(g\sharp f)_i$ by $H_i:M\times[s_{i-1},s_i]\rightarrow\mathbb{R}$, then
\begin{eqnarray*}
&&H_i(x,s)=2F_i(x,2s),\quad 1\leq i\leq n,\\
&&H_i(x,s)=2G_{i-n}(x,2s-1),\quad n+1\leq i\leq n+m.
\end{eqnarray*}
By definition we have
\begin{eqnarray*}
\begin{split}
&d_k(\phi,\varphi)\leq {\rm Length}_k(g\sharp f)
=\sum_{j=1}^{n+m}{\rm Length}_k\Big((g\sharp f)_j\Big) \\
=&\sum_{j=1}^{n+m}\sum_{i=0}^{k}\int_{s_{j-1}}^{s_j}\left\|\frac{\partial^i(H_j(x,s))}
{\partial s^i}\right\|_\infty ds \\
=&\sum_{j=1}^n\sum_{i=0}^{k}\int_{s_{j-1}}^{s_j}\left\|\frac{\partial^i(2F_j(x,2s))}
{\partial s^i}\right\|_\infty ds\\
& +
\sum_{j={n+1}}^{n+m}\sum_{i=0}^{k}\int_{s_{j-1}}^{s_j}\left\|\frac{\partial^i(2G_{j-n}
(x,2s-1))}{\partial s^i}\right\|_\infty ds
\end{split}
\end{eqnarray*}
\begin{eqnarray*}
\begin{split}
=&\sum_{j=1}^n\sum_{i=0}^{k}\int_{\frac{t_{j-1}}{2}}^{\frac{t_j}{2}}2^{i+1}\left\|
(F_j)_2^{(i)}(x,2s)\right\|_\infty ds \\
&\qquad+
\sum_{j={n+1}}^{n+m}\sum_{i=0}^{k}\int_{\frac{t_{j-1-n}'+1}{2}}^{\frac{t_{j-n}'+1}{2}}
2^{i+1}\left\|(G_{j-n})_2^{(i)}(x,2s-1)\right\|_\infty ds \\
=&\sum_{j=1}^n\sum_{i=0}^{k}\int_{t_{j-1}}^{t_j}2^i\left\|
(F_j)_2^{(i)}(x,s)\right\|_\infty ds + \sum_{j=1}^m\sum_{i=0}^{k}\int_{t_{j-1}'}^{t_j'}
2^i\left\|(G_j)_2^{(i)}(x,s)\right\|_\infty ds \\
\leq& 2^k\left(\sum_{j=1}^n\sum_{i=0}^{k}\int_{t_{j-1}}^{t_j}\left\|
\frac{\partial^i(F_j(x,s))}{\partial s^i}\right\|_\infty ds +
\sum_{j=1}^m\sum_{i=0}^{k}\int_{t_{j-1}'}^{t_j'}
\left\|\frac{\partial^i(G_j(x,s))}{\partial s^i}\right\|_\infty ds\right) \\
=&2^k\left(\sum_{j=1}^n{\rm Length}_k(f_j)+\sum_{j=1}^m{\rm Length}_k(g_j)\right) \\
=&2^k\left({\rm Length}_k(f)+{\rm Length}_k(g)\right).
\end{split}
\end{eqnarray*}
Take the infimum for all $f\in\Omega(\phi,\theta)$ and
$g\in\Omega(\theta,\varphi)$ at the right hand of the above equation respectively,
we obtain the desired triangle inequality $d_k(\phi,\varphi)\leq
2^k(d_k(\phi,\theta)+d_k(\theta,\varphi))$.\

\noindent{(iii)} Because $d_k\geq d_0=d_{\rm H}$, the non-degeneracy of $d_k$ is
obvious.\

\noindent{(iv)} Firstly, we prove the right-invariance of $d_k$. Let
$f\in\Omega(\phi,\varphi)$ be generated by $\{F_i,1\leq i\leq n\}$
as above. Then $(f\circ\theta)(t):=f(t)\theta$ defines an element
$f\circ\theta$ in $\Omega(\phi\theta, \varphi\theta)$, and
$f\circ\theta|_{[t_{i-1},t_i]}$ also correspond to (time-dependent)
Hamiltonian functions $F_i$, $i=1,\cdots,n$. So ${\rm
Length}_k(f)={\rm Length}_k(f\circ\theta)$. Hence,
\begin{eqnarray*}
d_k(\phi,\varphi) &=&\inf_f\{{\rm Length}_k(f)\}=
\inf_f\{{\rm Length}_k(f\circ\theta)\}\\
&=&\inf_{g\in\Omega(\phi\theta,\varphi\theta)} \{{\rm
Length}_k(g)\}=d_k(\phi\theta,\varphi\theta).
\end{eqnarray*}
 Next we prove the left-invariance of $d_k$. Let
$f$ as above. Then $\theta\circ f(t):=\theta(f_t(x))$ defines an
element $\theta\circ f$ in $\Omega(\theta\phi, \theta\varphi)$, and
$\theta\circ f|_{[t_{i-1},t_i]}$  correspond to Hamiltonian functions
$F_i(\theta^{-1}(\cdot),t)$, $i=1,\cdots,n$. So ${\rm
Length}_k(f)={\rm Length}_k(\theta\circ f)$. By the same argument as
above, we get $d_k(\phi,\varphi)=d_k(\theta\phi,\theta\varphi)$.
\hfill$\Box$\vspace{2mm}

\subsection{Proof of Theorems~\ref{th:HoferE},~\ref{th:SikE}}

Following \cite{Sik90, Ho93} we first prove

\begin{lemma}\label{lem:sikH}
Assume $\psi_1,\psi_2,\cdots,\psi_m\in{\rm Ham}(\mathbb{R}^{2n},\omega_0)$ have properly
separated supports. Then
$$
\|\psi_1\psi_2\cdots\psi_m\|^\ast_k\leq2(k+1)\max_j\|\psi_j\|_k^\ast.
$$
\end{lemma}

\noindent{\bf Proof.}
Given $\varepsilon>0$, by (\ref{def:k-dis*}) we have $f_j\in\Omega({\rm id},\psi_j)$ such that
${\rm Length}^\ast_k(f_j)<\|\psi_j\|_k^\ast+\varepsilon$, $j=1,\cdots,m$.
Let $H_j$ be the corresponding normalized Hamiltonian functions of $f_j$, $j=1,\cdots,m$.
 Then $|||H_j|||_k={\rm Length}^\ast_k(f_j)\leq \|\psi_j\|_k^\ast+\varepsilon$.
Through a refinement, we could suppose that there exists a common division of time
$0=t_0<t_1<\cdots<t_n=1$ such that $H_j^i:=H_j|_{[t_{i-1},t_i]}$ is smooth for every
$1\leq j\leq m,1\leq i\leq n$.

Using the same notations as those of \cite{HoZe}, we set
 $S_j={\rm supp}(\psi_j)$, choose $R>0$ such that ${\rm supp}(H_j)\subset B_R(x_j^\ast)\times
  [0, 1]$ with $x_j^\ast\in S_j$, and then vectors $v_j$ such that the sets $B_R(S_j+v_j)$
  are disjoint.  Let $\tau\in{\rm Ham}(M,\omega)$ be the map associated to the $S_j$
  guaranteed by  Lemma~8 on the page 175 of \cite{HoZe}, and let $\hat{\psi}_j=\tau\psi_j\tau^{-1}$ and
$$\hat{f}_j(t)=\tau f_j(t)\tau^{-1},\;j=1,\cdots,m.$$
 Then the corresponding normalized Hamiltonian function with $\hat{f}_j$ is given by
 $$
\hat{H}_j(x,t)=H_j(\tau^{-1}(x),t)=H_j(x-v_j,t),
$$
and thus  the corresponding normalized Hamiltonian function with
$$\hat{f}_1\hat{f}_2\cdots\hat{f}_m\in\Omega({\rm id}, \hat{\psi}_1\hat{\psi}_2\cdots\hat{\psi}_m)
$$
is $\hat{H}=\hat{H}_1+\cdots+\hat{H}_m$. It follows that
$$
\|\psi_1\psi_2\cdots\psi_m\|_k^\ast=\|\hat{\psi}_1\hat{\psi}_2\cdots\hat{\psi}_m\|_k^\ast
\leq|||\hat{H}|||_k.
$$
By the definition in (\ref{def:k-length*}) we have $s_{i}\in [0,1]$, $i=0,\cdots,k$,
such that
\begin{eqnarray*}
|||\hat{H}|||_k&=&\sum_{i=0}^k\max_{1\le l\le n}\max_{t_{l-1}\le t\le t_l}\left(\sup_x\frac{\partial^i\hat{H}^l}{\partial t^i}
-\inf_x\frac{\partial^i\hat{H}^l}{\partial t^i}\right)\\
               &=&\sum_{i=0}^k\left(\sup_x\frac{\partial^i\hat{H}}{\partial t^i}(x,s_i)
-\inf_x\frac{\partial^i\hat{H}}{\partial t^i}(x,s_i)\right)\\
               &=&\sum_{i=0}^k\left[\max_{1\leq j\leq m}\left(\sup_x\frac{\partial^i\hat{H}_j
}{\partial t^i}(x,s_i)\right)-\min_{1\leq j\leq m}\left(\inf_x\frac{\partial^i\hat{H}_j}
{\partial t^i}(x,s_i)\right)\right]\\
               &\leq&\sum_{i=0}^k2\max_{1\leq j\leq m}\left(\left\|\frac{\partial^i\hat{H}_j}
{\partial t^i}(\cdot,s_i)\right\|_\infty\right)\\
               &\leq&\sum_{i=0}^k2\max_{1\leq j\leq m}\left(\max_{t\in [0,1]}\left\|\frac{\partial^i\hat{H}_j
}{\partial t^i}(\cdot,t)\right\|_\infty\right)\\
               &\leq&2(k+1)\max_{1\leq j\leq m}\left(\sum_{i=0}^k\max_{t\in [0,1]}\left\|
\frac{\partial^i\hat{H}_j}{\partial t^i}(\cdot,t)\right\|_\infty\right)\\
               &=&2(k+1)\max_{1\leq j\leq m}|||\hat{H}_j|||_k\\
               &\leq&2(k+1)\max_j\|\psi_j\|_k^\ast+2(k+1)\varepsilon.
\end{eqnarray*}
This holds for every $\varepsilon>0$ and the lemma is proved.
\hfill$\Box$\vspace{2mm}

\noindent{\bf Proof of Theorem~\ref{th:SikE}}.
Given $\varepsilon>0$, as in the proof of \cite{Sik90, Ho93} we can construct maps $\psi_j,
0\leq j\leq N$ satisfying
$d_k^\ast(\psi_j,\psi_{j+1})<\varepsilon$,
and maps $\varphi_j(0\leq j\leq2N)$, $\alpha_j(1\leq j\leq N)$, $\beta_j(0\leq j\leq N)$ with
$\|\varphi_0\|_k^\ast\leq e_{p,k}^\ast(U)+\varepsilon$.
Then we have
\begin{eqnarray*}
\|\varphi_H\|_k^\ast&=&\|\beta_N\|_k^\ast=\left\|(\prod_{j=1}^N\alpha_j\beta_j)
(\prod_{j=0}^{N-1}\alpha_{j+1}\beta_j)^{-1}\right\|_k^\ast\\
                &\leq&2^{k+2}(k+1)\left(\max_{1\leq j\leq N}\|\alpha_j\beta_j\|_k^\ast+
\max_{0\leq j\leq N-1}\|\alpha_{j+1}\beta_j\|_k^\ast\right)
\end{eqnarray*}
(because $\|\phi\psi\|^\ast_k\le 2^{k+1}(\|\phi\|^\ast_k+\|\psi\|^\ast_k)$ by Theorem~\ref{th:corres}).
We can estimate
\begin{eqnarray*}
\|\alpha_j\beta_j\|_k^\ast&\leq&2^{k+1}\left[\|\varphi_{2j-1}\|_k^\ast+2^{k+1}(
\|\varphi_{2j-1}^{-1}\varphi_{2j}\|^\ast_k+\|\varphi_{2j}\|_k^\ast)\right]\\
                      &\leq&2^{k+1}\|\varphi_0\|_k^\ast+2^{3k+3}(\|\varphi_{2j-1}\|_k^\ast+
\|\varphi_{2j}\|_k^\ast)+2^{2k+2}\|\varphi_0\|_k^\ast\\
                      &=&(2^{k+1}+2^{2k+2}+2^{3k+4})\|\varphi_0\|_k^\ast.
\end{eqnarray*}
Similarly we have
\begin{eqnarray*}
\|\alpha_{j+1}\beta_j\|_k^\ast&\leq&2^{k+1}\left[\|\varphi_{2j+1}\|_k^\ast+2^{k+1}
(\|\psi_{j+1}\psi_j^{-1}\psi_j\varphi_{2j+1}^{-1}\varphi_{2j}\psi_j^{-1}\|_k^\ast
+\|\varphi_{2j}\|_k^\ast)\right]\\
                          &\leq&2^{k+1}\|\varphi_0\|_k^\ast+2^{3k+3}(\|\psi_{j+1}\psi_j^{-1}\|_k^\ast
+\|\varphi_{2j+1}^{-1}\varphi_{2j}\|_k^\ast)+2^{2k+2}\|\varphi_0\|_k^\ast\\
                          &\leq&(2^{k+1}+2^{2k+2}+2^{4k+5})\|\varphi_0\|_k^\ast+2^{3k+3}
d_k^\ast(\psi_{j+1},\psi_j).
\end{eqnarray*}
Summing up, we have
\begin{eqnarray*}
\|\varphi_H\|_k^\ast&\leq&2^{k+2}(k+1)(2^{k+2}+2^{2k+3}+2^{3k+4}+2^{4k+5})\|\varphi_0\|_k^\ast+2^{4k+5}
(k+1)\varepsilon\\
                 &\leq&2^{2k+4}(k+1)(1+2^{k+1}+2^{2k+2}+2^{3k+3})e_{p,k}^\ast(U)+\\
&&2^{2k+4}(k+1)(1+2^{k+1}+2^{2k+2}+2^{3k+3})\varepsilon+2^{4k+5}(k+1)\varepsilon.
\end{eqnarray*}
Since $\varepsilon>0$ is arbitrary, the proof is finished.
\hfill$\Box$\vspace{2mm}

Carefully checking the proof of Proposition 6 in \cite{Ho93} (or \cite[Lemma 10]{HoZe})
and replacing $E$ and $e_p$ therein
we still have for our $\|\cdot\|_k^\ast$ and $e_{p,k}^\ast$:

\begin{lemma}\label{lem:L-Sik}
Let $\psi\in{\rm Ham}(\mathbb{R}^{2n},\omega_0)$ with $\psi\neq{\rm id}$ and let
$\delta>|\psi-{\rm id}|_{C^0}$. For every $Q\subset\mathbb{R}^{2n}$ open and satisfying
$Q\cap{\rm supp}(\psi)\neq\emptyset$, there exists a $\varphi\in{\rm Ham}(\mathbb{R}^{2n},
\omega_0)$ satisfying
\begin{description}
\item[(i)] $\varphi|Q=\psi|Q$
\item[(ii)] ${\rm supp}(\psi)\subset U$
\item[(iii)] $\|\varphi\|_k^\ast\leq2^{2k+4}(k+1)(1+2^{k+1}+2^{2k+2}+2^{3k+3})
e_{p,k}^\ast(U)$,
\end{description}
where $U$ is the intersection of $B_\delta(Q)$ with the convex hull of ${\rm supp}(\psi)$, and
$B_\delta(Q)=\{x|{\rm dist}(x,Q)<\delta\}$.
\end{lemma}

Similarly, corresponding to \cite[Corollary 7]{Ho93} or \cite[Lemma 11]{HoZe} we have

\begin{lemma}\label{lem:H}
Let $U=(a_1,a_2)\times(b_1,b_2)\oplus\mathbb{R}^{2n-2}$.
Then for $\psi\in{\rm Ham}(\mathbb{R}^{2n},\omega_0)$ with ${\rm supp}(\psi)\subset U$ it holds that
$$
\|\varphi\|_k^\ast\leq2^{2k+4}(k+1)(1+2^{k+1}+2^{2k+2}+2^{3k+3})(a_2-a_1)(b_2-b_1).
$$
\end{lemma}

\noindent{\bf Proof of Theorem~\ref{th:HoferE}}.
Let $\{\varphi_j\}_{j\in I}$ be as in the proof of \cite[Theorem 9]{HoZe}. Corresponding
to the inequality (iii) on the page 180 of \cite{HoZe}, we  have
$$
\|\varphi_j\|_k^\ast\leq2^{2k+4}(k+1)(1+2^{k+1}+2^{2k+2}+2^{3k+3})e_{p,k}^\ast(U_j).
$$
And similar to \cite[(5.32)]{HoZe} we have estimate
\begin{eqnarray*}
\|\prod_{j\in I}\varphi_j\|_k^\ast&\leq&2(k+1)\max_{j\in I}\|\varphi_j\|_k^\ast\\
&\leq&2(k+1)2^{2k+4}(k+1)(1+2^{k+1}+2^{2k+2}+2^{3k+3})\cdot2(\varepsilon+\delta_1)R\\
&\leq&2^{2k+6}(k+1)^2(1+2^{k+1}+2^{2k+2}+2^{3k+3})R\delta.
\end{eqnarray*}
Thus
\begin{eqnarray*}
\|\psi\|_k^\ast&\leq&2^{k+1}(\|\psi\theta\|_k^\ast+\|\theta\|_k^\ast)\\
&=&2^{k+1}(\|\prod\hat{\varphi}_j\|_k^\ast+\|\prod\varphi_j\|_k^\ast)\\
&\leq&2^{3k+8}(k+1)^2(1+2^{k+1}+2^{2k+2}+2^{3k+3})R\delta.
\end{eqnarray*}
This holds true for every $\delta>|\psi-{\rm id}|_{C^0}$, and we conclude that
$$
\|\psi\|_k^\ast\leq2^{3k+8}(k+1)^2(1+2^{k+1}+2^{2k+2}+2^{3k+3})R|\psi-{\rm id}|_{C^0},
$$
where $R={\rm diameter}\bigl({\rm supp}(\psi)\bigr)$. Finally we have the inequality:
$$
d_k^\ast(\varphi,\psi)\leq2^{3k+8}(k+1)^2(1+2^{k+1}+2^{2k+2}+2^{3k+3})\,{\rm diameter}\bigl({\rm supp}(\varphi\psi^{-1})\bigr)
|\varphi-\psi|_{C^0},
$$
\hfill$\Box$\vspace{2mm}

\subsection{Proof of Theorem~\ref{th:EP}}

When $k=0$ and $p\in [1,\infty)$ this is the result in \cite{EP93}. We shall assume
that either $k=0$ and $0<p<1$ or $k\in\mathbb{N}$ and $0<p<1/k$ below.

The proof ideas are same as those of \cite{EP93}.
Given  a bi-invariant pseudo quasimetric $\rho$ on ${\rm Ham}(M, \omega)$,
 the ($\rho$-)\textbf{displacement energy} of a subset $A\subset M$  is defined by
\begin{equation}\label{e:disE}
e_\rho(A)=\inf\{\rho({\rm id},g)|\;g\in {\rm Ham}(M, \omega)\;\hbox{such
that}\; g(A)\cap A=\emptyset\}.
\end{equation}
As in the proof of \cite[Th.2.2.A]{EP93} we may obtain:

\begin{claim}\label{cl:postiveDE}
Every nonempty open subset of $M$ has positive displacement energy with respect to
 a bi-invariant quasimetric on ${\rm Ham}(M, \omega)$.
\end{claim}

So it suffices to prove that  the
displacement energy associated with $d_{(k,p)}$ ($0<p<1/k$) vanishes
for some embedded open ball in $M$. In fact,
using Darboux theorem, we can choose a chart $M\supset U\ni w\to (x_1(w),\cdots,x_n(w),y_1(w),\cdots, y_n(w))\in\mathbb{R}^{2n}$ so that
the symplectic form $\omega=\sum_{i=1}^ndx_i\wedge
dy_i$ on it. Replacing $\omega$ by $N\omega$ for some large $N>0$ we may assume
\begin{eqnarray*}
U\supset B(0,4)=\left\{(x,y)\in\R^{2n}\,\Bigm|\, \sum^n_{i=1}(x_i^2+y_i^2)<
16\right\},\quad A=B(0,1).
\end{eqnarray*}

Consider a Hamiltonian isotopy $\{h_t\}$, $t\in[0,1]$ such that when
restricted to $U$, $h_t$ is simply a shift by $2t$ along the $y_1$
coordinate. Assume that $\{h_t\}$ is generated by $H$, then
$H(x,y,t)=2x_1$ on $U$ ($H$ is not normalized). Clearly $h_1(A)\cap
A=\emptyset$.

Fixed a smooth cut-off function $\delta:\mathbb{R}\rightarrow[0,1]$,
such that $\delta\equiv1$ for $|t|\leq1/4$, $\delta\equiv0$ for
$|t|\geq3/4$. Set $\delta_m(t)=\delta(mt)$. When $M$ is open, we can
define a sequence of functions $\{G_m\}_{m\in\N}$ as follows:
\begin{equation}\label{e:G_m}
G_m(x,y,t)=
\begin{cases}
2x_1\delta_m(\sqrt{|x|^2+|y-(2t,0,\cdots,0)|^2}-1) & (x,y)\in U\\
0 & (x,y)\notin U.
\end{cases}
\end{equation}
$G_m$ is smooth on $M$. The support of $G_m$ is contained in a
neighborhood of $h_t(\partial A)$ for each $t\in[0,1]$ and tends to
$h_t(\partial A)$ as $m\rightarrow\infty$. Since for every $t$ the
function $G_m(\cdot,t)$ coincides with $H$ near $h_t(\partial A)$,
we conclude that the Hamiltonian vector field of $G_m$ is equal to
$X_H$ near $h_t(\partial A)$ for every $m$. Hence the Hamilton
isotopy $\{\phi_{G_m}^t\}$ of $G_m$ satisfies $\phi_{G_m}^t(\partial
A)=h_t(\partial A)$ and so $\phi_{G_m}(\partial A)\cap\partial
A=\emptyset$. But this obviously implies that $\phi_{G_m}(A)\cap
A=\emptyset$.

Clearly, for $0<p<\infty$ we have $\int^1_0\|G_m(\cdot,t)\|_pdt\to 0$ as $m\to\infty$.

Next we prove that  for each $1\leq i\leq k$,
$$\left\|\frac{\partial^iG_m}{\partial
t^i}(\cdot,t)\right\|^p_p \rightarrow
0\quad\hbox{uniformly in $t\in [0,1]$ as}\;m\rightarrow\infty
$$
 For simplicity we write
$r=\sqrt{(|x|^2+|y-(2t,0,\cdots,0)|^2)}$. In $U$ we have
$$
\frac{\partial^iG_m}{\partial
t^i}(x,y,t)=m^ix_1\delta_m^{(i)}(r-1)\cdot
2^{i+1}\big(\frac{t-y_1}{r}\big)^i+m^{i-1}L(x,y,t).
$$
Here $L$ is a combination of $\delta_m^{(j)}, 1\leq j\leq i-1$. Its
coefficients are rational functions of $x,y,t$ which are bounded in
$C^0$-norm near $h_t(\partial A)$. Hence,
\begin{eqnarray*}
\Big\|\frac{\partial^iG_m}{\partial t^i}(\cdot,t)
\Big\|_p^p&=&\int_Um^{ip}\Big|x_1\delta_m^{(i)}(r-1)\cdot
2^{i+1}\left(\frac{t-y_1}{r}\right)^i+\frac{1}{m}L(x,y,t)\Big|^p\omega^n\\
&\leq& Cm^{ip}\int_{\Sigma_1\cup \Sigma_2}\omega^n
=Cm^{ip}\bigl({\rm Vol}(\Sigma_1)+{\rm Vol}(\Sigma_2)\bigr),
\end{eqnarray*}
where $C>0$ is a constant depending on $\delta$, $\Sigma_1$ is a
region bounded by two spheres whose radiuses are $1+\frac{1}{4m},
1+\frac{3}{4m}$ respectively, and $\Sigma_2$ is a similar region
with radiuses $1-\frac{1}{4m},1-\frac{3}{4m}$ respectively. So we
have
\begin{eqnarray*}
&&\left\|\frac{\partial^iG_m}{\partial t^i}(\cdot,t)
\right\|_p^p\\
&\leq&\frac{\pi^n}{n!}Cm^{ip}
\left[\left(1+\frac{3}{4m}\right)^{2n}-\left(1+\frac{1}{4m}\right)^{2n}+
\left(1-\frac{1}{4m}\right)^{2n}-\left(1-\frac{3}{4m}\right)^{2n}\right]\\
&=&\frac{\pi^n}{n!}Cm^{ip}\left[\frac{2}{m}+O\left(\frac{1}{m^2}\right)\right].
\end{eqnarray*}
Since $ip\leq kp<1$, so the above expression tends to zero when
$m\rightarrow\infty$.

When $M$ is closed, we could find an open set $V$ disjoint with the
above $U$, and shrinking $U$ properly we can also assume that $V$ is
symplectomorphic to $U$. Then we can define a function
$\widetilde{G_m}$ on $V$ which is the same form as $G_m$ in
(\ref{e:G_m}) but with a minus in front additionally. Define $K_m$
to be the sum of $G_m$ and $\widetilde{G_m}$. Then $K_m$ is a
normalized Hamiltonian function such that $\phi_{K_m}(A)\cap
A=\emptyset$. Using the same estimates as above for $G_m$ and
$\widetilde{G_m}$ respectively, we derive
$$
\left\|\frac{\partial^iK_m}{\partial t^i}(\cdot,t)\right\|_p
\rightarrow 0\quad\hbox{uniformly in $t\in [0,1]$ as}\quad m\rightarrow\infty.
$$
All these lead to the first conclusion.

Note that $\{\phi\in {\rm Ham}(M, \omega)\,|\, \rho({\rm id}_M,\phi)=0\}$
is also a normal subgroup of ${\rm Ham}(M, \omega)$
for any bi-invariant pseudo quasimetric $\rho$ on ${\rm Ham}(M, \omega)$.
The second claim follows from Banyaga's theorem as in \cite{PoBook}.

\hfill$\Box$\vspace{2mm}

\section{Extensions}\setcounter{equation}{0}

It is a natural question to extend bi-invariant
pseudo-metric on ${\rm Ham}(M,\omega)$  to the group of symplectic
diffeomorphisms. There exist different ways realizing this. We only use the method by
Lalonde and Polterovich \cite{LaPo97}, and one by  Banyaga \cite{Ban07}.

\subsection{The method by Lalonde and Polterovich}

 Each  $\phi\in{\rm Symp}(M,\omega)$ induces an
isometry with respect to the Hofer metric,
$$
C_\phi:{\rm Ham}(M,\omega)\rightarrow {\rm Ham}(M,\omega),\;f\mapsto
C_\phi f:=\phi f\phi^{-1}.
$$
For $\alpha\in (0, \infty]$ it was shown in \cite[Proposition
1.2.A]{LaPo97} that
\begin{equation}\label{e:alpha-bouned}
r_\alpha(\phi):=\sup\{d_H(f,C_\phi f)\,|\, f\in{\rm
Ham}(M,\omega)\;\&\;\|f\|_H\le\alpha\}
\end{equation}
defines a bi-invariant function $r_\alpha$ and $r_\alpha(\phi)\le
2\alpha\;\forall\phi\in {\rm Symp}(M,\omega)$. In particular,
$r_\alpha$ is an \textsf{bi-invariant norm} on ${\rm
Symp}(M,\omega)$ if $\alpha\in (0,\infty)$. We say
$\phi\in{\rm Symp}(M,\omega)$ to be {\bf bounded} if the function
$(0,\infty)\ni\alpha\mapsto r_\alpha(\phi)$ is bounded, or equivaliently
$C_\phi$ is $C^0$-bounded, i.e.,
\begin{equation}\label{e:bouned}
r_\infty(\phi)=\sup\{d_H(f,C_\phi f)\,|\, f\in{\rm
Ham}(M,\omega)\}<\infty.
\end{equation}
By the bi-invariance of the Hofer metric $d_H$ the set of all
bounded symplectomorphism  in ${\rm Symp}(M,\omega)$ form  a normal
subgroup of ${\rm Symp}(M,\omega)$, denoted by ${\rm BI}(M)$. Set
${\rm BI}_0(M)={\rm BI}(M)\cap \rm{Symp}_0(M)$, where
$\rm{Symp}_0(M)$ is the connected component of ${\rm
Symp}(M,\omega)$ containing the identity map.
 Every
$\phi\in {\rm Ham}(M,\omega)$ is bounded and $r_\infty(\phi)\le
2\|\phi\|_H$ because
\begin{equation}\label{e:b-i ineq}
d_H(f,C_\phi f)=d_H({\rm id},[\phi,f])\leq
2\|\phi\|_H,\quad\forall f\in{\rm Ham}(M,\omega).
\end{equation}
Hence ${\rm Ham}(M,\omega)$ is a subgroup of ${\rm BI}_0(M)$.
These motivated Lalonde and Polterovich  to propose the so-called
 \textsf{bounded isometry conjecture}.

\begin{conjecture}[\hbox{\cite[Conjecture
1.3.A]{LaPo97}}]\label{conj:LP}
 ${\rm BI}_0(M)={\rm Ham}(M,\omega)$.
\end{conjecture}

They proved it in \cite{LaPo97} for some cases, for example, $M$ is
any closed surface with area form or $M$ is a product of closed
surfaces of genus greater than $0$ with product symplectic form.
For recent progresses the reader may refer to \cite{LaPe, Han, CaPe,Pedroza}.

Now let us consider corresponding questions with metric $d_k$.
 For each
$\phi\in{\rm Symp}(M,\omega)$ it is easy to check that $C_\phi$ is
still an isometry of ${\rm Ham}(M,\omega)$ with respect to the quasi-norm
 $\|\cdot\|_k= d_k({\rm id},\cdot)$.
Corresponding to (\ref{e:alpha-bouned}) we define
\begin{equation}\label{e:k-alpha-bouned}
r_{\alpha,k}(\phi):=\sup\{d_k(f,C_\phi f)\,|\, f\in{\rm
Ham}(M,\omega)\;\&\;\|f\|_k\le\alpha\}
\end{equation}
for each $\alpha\in (0, \infty]$. It is also bi-invariant and
satisfies
\begin{equation}\label{e:k-alpha-bound}
r_{\alpha,k}(\phi)\le 2^{k+1}\alpha,\;\forall\phi\in{\rm Symp}(M,\omega).
\end{equation}
Actually,  the function $r_{\alpha,k}$ is a quasi-norm on ${\rm Symp}(M,\omega)$
by the following:

\begin{lemma}\label{lem:biInvariance}
For every $k\geq0$, $\alpha\in(0,\infty]$, the function $r_{\alpha,k}$ on ${\rm
Symp}(M,\omega)$ defined above is conjugate invariant, assumes the
value $0$ only at the identity, and satisfies the quasi-triangle
inequality
$$
r_{\alpha,k}(\phi\psi)\leq 2^k(r_{\alpha,k}(\phi)+r_{\alpha,k}(\psi)).
$$
\end{lemma}

\noindent{\bf Proof.} (i) \underline{The conjugate invariance of
$r_{\alpha,k}$}. For $\phi,\varphi\in{\rm Symp}(M,\omega)$ we have
\begin{eqnarray*}
r_{\alpha,k}(\varphi\phi\varphi^{-1})&=&\sup_{\{f|\|f\|_k\leq\alpha\}}d_k
(f,\varphi\phi\varphi^{-1}f\varphi\phi^{-1}\varphi^{-1})\\
&=&\sup_{\{f|\|f\|_k\leq\alpha\}}
d_k(\varphi^{-1}f\varphi,\phi\varphi^{-1}f\varphi\phi^{-1})\\
&=&\sup_{\{g|\|g\|_k\leq\alpha\}}d_k(g,\phi
g\phi^{-1})\qquad\quad\hbox{(setting $g=\varphi^{-1}f\varphi$)}\\
&=&r_{\alpha,k}(\phi).
\end{eqnarray*}

\noindent{(ii)} \underline{The non-degeneracy of $r_{\alpha,k}$}. If
$r_{\alpha,k}(\phi)=0$, then $d_k(f,\phi f\phi^{-1})=0$ for each $f\in{\rm
Ham}(M,\omega)$, $\|f\|_k\leq\alpha$. By the non-degeneracy of $d_k$, we have $f=\phi
f\phi^{-1}$. Suppose that $f$ is generated by Hamiltonian function
$F$. Then $\phi f\phi^{-1}$ is generated by $F\circ\phi^{-1}$. Hence
$F=F\circ\phi^{-1}$. By the arbitrariness of $F$, we have
$\phi={\rm id}$.

\noindent{(iii)} \underline{The quasi-triangle inequality holds for
$r_{\alpha,k}$}.
\begin{eqnarray*}
r_{\alpha,k}(\phi\psi)&=&\sup_{\{f|\|f\|_k\leq\alpha\}}d_k(f,\phi\psi
f\psi^{-1}\phi^{-1})\\
&\leq& 2^k\left(\sup_{\{f|\|f\|_k\leq\alpha\}}d_k(f,\phi
f\phi^{-1})+
\sup_{\{f|\|f\|_k\leq\alpha\}}d_k(\phi f\phi^{-1},\phi\psi f\psi^{-1}\phi^{-1})\right)\\
&=& 2^k\left(r_{\alpha,k}(\phi)+r_{\alpha,k}(\psi)\right).
\end{eqnarray*}
\hfill$\Box$\vspace{2mm}

If $r_{\infty, k}(\phi)<\infty$  we call $\phi$ a \textbf{k-bounded
isometry}. Clearly, every $k$-bounded symplectomorphism is bounded.
As in (\ref{e:b-i ineq}) for $\phi\in {\rm Ham}(M,\omega)$ we have
\begin{equation}\label{e:k-b-i ineq}
d_k(f,C_\phi f)=d_k({\rm id},[\phi,f])\leq
2^{k+1}\|\phi\|_k,\quad\forall f\in{\rm Ham}(M,\omega)
\end{equation}
and thus  $r_{\infty,k}(\phi)\le 2^{k+1}\|\phi\|_k$.

Let
\begin{eqnarray*}
&&{\rm BI}_k(M)=\{\phi\in {\rm Symp}(M,\omega)\,|\,
r_{\infty,k}(\phi)<\infty\}\quad\hbox{and}\quad\\
&& {\rm BI}_{k,0}(M)={\rm BI}_k(M)\cap \rm{Symp}_0(M),
\end{eqnarray*}
which correspond to the sets ${\rm BI}(M)$ and ${\rm BI}_{0}(M)$
respectively. It is obvious that
$$
{\rm BI}(M)\supseteq {\rm BI}_1(M)\supseteq{\rm
BI}_2(M)\supseteq\cdots\supseteq{\rm
BI}_k(M)\supseteq\cdots\supseteq {\rm Ham}(M,\omega)
$$
because (\ref{e:k-b-i ineq}) and
$\|\cdot\|_H=\|\cdot\|_0\le\|\cdot\|_1\le\cdots\le\|\cdot\|_k\le\cdots$.
It is probably that ${\rm BI}_k(M)\ne{\rm BI}(M)$ for some $k\in\N$
and $M$. Corresponding to Conjecture~\ref{conj:LP} we have

\begin{conjecture}\label{conj:weak-LP}
For every symplectic manifold $M$ and every integer $m\in\N\cup\{0\}$,
$$
\bigcap_{k=m}^{\infty}{\rm BI}_{k,0}(M)={\rm
BI}_{m-1,0}(M).\eqno({\rm BIC})_m
$$
In particular, for every symplectic manifold $M$,
$$
\bigcap_{k=0}^{\infty}{\rm BI}_{k,0}(M)={\rm
Ham}(M,\omega).\eqno({\rm WBIC})
$$
\end{conjecture}
Here (WBIC) means \textsf{weak bounded isometry conjecture}.
Clearly, the proof of ({\rm WBIC}) is more hopeful  than one of
Conjecture~\ref{conj:LP}.  In particular, all manifolds mentioned
above satisfy $({\rm BIC})_m$ and $({\rm WBIC})$.

In the following we shall point out that many results in
\cite{LaPo97} can be generalized to the case of our quasi-metrics.
Carefully checking the proof of Theorem~1.3.F in \cite{LaPo97}  we
immediately obtain the following generalization of it.

\begin{proposition}\label{prop:1.3.F}
Suppose that a symplectic manifold $(M,\omega)$ satisfies ${\rm
BI}_{k,0}(M)={\rm Ham}(M,\omega)$ for some $k\in\N\cup\{0\}$. Then
every symplectomorphism $\phi\in{\rm Symp}(M,\omega)$ which acts
nontrivially on $H_c^1(M, \R)$ is $k$-unbounded.
\end{proposition}

Since a $k$-bounded symplectomorphism is also bounded, the following
two propositions are, respectively, direct consequences of
Theorems~1.4.A,~1.3.C in \cite{LaPo97}.

\begin{proposition}\label{prop:1.4.A}
If a closed Lagrangian submanifold $L\subset (M, \omega)$
admits a Riemannian metric with non-positive sectional curvature, and
 the inclusion $L\hookrightarrow M$ induces an injection
on fundamental groups, then  $\phi(L)\cap L\neq\emptyset$ for every
$k$-bounded $\phi\in{\rm Symp}(M,\omega)$.
\end{proposition}

\begin{proposition}\label{prop:1.3.C}
 Let $S$ be a closed surface of genus greater than $0$ and let $(W,\omega_W)$ be
closed and weakly exact (i.e., $\omega_W|_{\pi_2(W)}=0$).
Suppose that $\phi\times\psi$ be a $k$-bounded symplectomorphism of
 $(S\times W, \omega_S\oplus\omega_W)$. Then the
symplectomorphism $\phi$ is Hamiltonian.
\end{proposition}

Finally, we give the corresponding result of Theorem~5.1.A in
\cite{LaPo97}.

\begin{thm}\label{th:BI(R^2n)}
For the standard symplectic space $(\mathbb{R}^{2n},\omega_0)$, any
compactly supported symplectic diffeomorphism  $\phi$ of
$\mathbb{R}^{2n}$ is $k$-bounded. Precisely, we have
$$
r_{\infty,k}(\phi)=\sup\{d_k(f,C_\phi f)\,|\,f\in{\rm Ham}(\mathbb{R}^{2n},\omega_0)\}\leq2^{3k+2}E({\rm supp}(\phi)),
$$
where $E({\rm supp}(\phi))$ is the cylindrical capacity of ${\rm supp}(\phi)$. Recall that cylindrical capacity $E(X)$ of a bounded
subset $X\subset\mathbb{R}^{2n}$ is definded by
$$
\inf\left\{c>0\,|\, \exists\;\phi\in {\rm Ham}(\R^{2n},\omega_0)\;
\hbox{such that}\;\phi(X)\subset
B^2(c)\times\mathbb{R}^{2n-2}\right\},
$$
where $B^2(c)$ stands for a disc of area $c$.
\end{thm}

This result can be proved  along the proof lines of
\cite[Theorem~5.1.A]{LaPo97}. For the sake of completeness we give
the proof of it. Firstly, we give a corresponding lemma  with
Lemma~5.1.B in \cite{LaPo97}.

\begin{lemma}\label{lem:estimate}
For all maps $\phi\in{\rm Symp}(\mathbb{R}^{2n},\omega_0)$ and $f,g\in{\rm
Ham}(\mathbb{R}^{2n},\omega_0)$,
$$
d_k(f,C_\phi f)\leq 2^kd_k(f,C_{g\phi g^{-1}}f)+2^{3k+2}\|g\|_k.
$$
\end{lemma}

\noindent{\bf Proof.}
The quasi-triangle inequality implies that
$$
d_k(f,C_\phi f)\leq 2^k(d_k(f,C_{g\phi g^{-1}}f)+d_k(C_\phi f,C_{g\phi g^{-1}}f)).
$$
Write $I=d_k(C_\phi f,C_{g\phi g^{-1}}f)$. Then
\begin{eqnarray*}
I&=&d_k(\phi f\phi^{-1},g\phi g^{-1}fg\phi^{-1}g^{-1})=d_k(f,\phi^{-1}g\phi
g^{-1}fg\phi^{-1}g^{-1}\phi)\\
&=&d_k(C_{[g,\phi^{-1}]}f,f),
\end{eqnarray*}
By Lemma~\ref{lem:biInvariance} we have
\begin{eqnarray*}
I&\leq&r_{\infty,k}([g,\phi^{-1}])=r_{\infty,k}(g\phi^{-1} g^{-1}\phi)\\
&\leq&2^k(r_{\infty,k}(g)+r_{\infty,k}(\phi^{-1}g^{-1}\phi))\\
&=&2^{k+1}r_{\infty,k}(g)\le 4^{k+1}\|g\|_k
\end{eqnarray*}
because of (\ref{e:k-b-i ineq}). \hfill$\Box$\vspace{2mm}

\begin{lemma}\label{lem:Sikorav}
Suppose that $K$ is a compact subset in $\mathbb{R}^{2n}$, and that
$L\subset\mathbb{R}^{2n}$ is a hyperplane so that $K$ lies on the
left of $L$. Let $g$ be a compactly supported Hamiltonian
diffeomorphism such that $g(K)$ sits in the right side of $L$, and
let $L'$ be  an arbitrary hyperplane parallel to $L$ such that
$g(K)$ lies between $L$ and $L'$. Then there exists another
compactly supported Hamiltonian diffeomorphism $g'$  such that
$\|g'\|_k=\|g\|_k$ and $g'(K)$ lies on the right of $L'$.
\end{lemma}

\noindent{\bf Proof.} Let $(x_1,y_1,\cdots,x_n,y_n)$ denote the
coordinates in $\R^{2n}$. Without loss of generality we may assume
that for some $v>0$ and $\varepsilon>0$,
\begin{eqnarray*}
&&L=\{x_1=0\}\quad\hbox{and}\quad L'=\{x_1=v\},\\
&&x_1<-\varepsilon\quad\forall (x_1,y_1,\cdots,x_n,y_n)\in K.
\end{eqnarray*}
Choose a cut off
function $\eta(t):\mathbb{R}\rightarrow[0,1]$ satisfying
$$
\eta|_{(-\infty,-\varepsilon]}=0\quad\hbox{and}\quad\eta|_{[0,+\infty)}=1.
$$
Let $S:\mathbb{R}^{2n} \rightarrow\mathbb{R}^{2n}$ be the
Hamiltonian diffeomorphism generated by the  function
$$
\R^{2n}\ni (x_1,y_1,\cdots,x_n,y_n)\mapsto H(x_1,y_1,\cdots,x_n,y_n)
=\eta(x_1)\cdot vy_1\in\R.
$$
(It is not compactly supported!) It is easily checked that
\begin{eqnarray*}
&&S((x_1,y_1,\cdots,x_n,y_n)=(x_1,y_1,\cdots,x_n,y_n)\quad\hbox{if}\;x_1<-\varepsilon,\\
&&S((x_1,y_1,\cdots,x_n,y_n)=(x_1+v,y_1,\cdots,x_n,y_n)\quad\hbox{if}\;x_1>0.
\end{eqnarray*}
Set $g':=SgS^{-1}$. It is also a compactly supported
Hamiltonian diffeomorphism and $\|g'\|_k=\|g\|_k$ by the
bi-invariance of the quasi-metric $d_k$. Clearly
$$
g'(K)=Sg
(K)=g(K)+(v,0,\cdots,0)
$$
lies on the right of $L'$. \hfill$\Box$\vspace{2mm}


\noindent{\bf Proof of Theorem~\ref{th:BI(R^2n)}.}  By the
definition of the cylindrical capacity, for a sufficiently small
$\delta>0$, there exists a Hamiltonian diffeomorphism $\psi$ such
that
$$
\psi({\rm supp}(\phi))\subset
B^2(c)\times\mathbb{R}^{2n-2},\quad\hbox{where}\quad c<E({\rm supp}(\phi))+\delta.
$$
By composing with a suitable Hamiltonian diffeomorphism we may
assume that  $\psi({\rm supp}(\phi))$ sits in
$Q\times\mathbb{R}^{2n-2}$, where $Q$ is an open square in the
$(x_1,y_1)$-plane with area $c$. Note that the displacement energy
of $Q$ is just $c$. By the example shown in \cite[p.17]{PoBook},
there exists a Hamiltonian isotopy $\{g_t\}\subset{\rm
Ham}(\R^2,\omega_0)$ with the time-$1$ map $g$ satisfying ${\rm
Length}(g_t)=c$ and $Q\cap g(Q)=\emptyset$. Since the Hamiltonian of
the flow $\{g_t\}$ may be chosen to be autonomous, we have
$\|g\|_k\leq{\rm Length}_k(g_t)={\rm Length}(g_t)=c$. Then for
$\widetilde{g}:=g\times id:\mathbb{R}^{2n}\rightarrow
\mathbb{R}^{2n}$, we have
$$
\|\widetilde{g}\|_k\leq\|g\|_k\leq
c<E({\rm supp}(\phi))+\delta.
$$

Obverse that  $\psi({\rm supp}(\phi))$ and $\widetilde{g}\psi({\rm supp}(\phi))$ have a positive distance. We can construct a
hyperplane $L$  lying between $\psi({\rm supp}(\phi))$ and
$\widetilde{g}\psi({\rm supp}(\phi))$ such that $\psi({\rm supp}(\phi))$ strictly sits in the left side of $L$. For an
arbitrary fixed $f\in{\rm Ham}(\mathbb{R}^{2n},\omega_0)$, since
${\rm supp}(f)$ is a compact set of $\mathbb{R}^{2n}$ by the
assumption we may choose a hyperplane $L'$
 parallel to $L$ such that $\widetilde{g}\psi({\rm supp}(\phi))$
 strictly lies between $L$ and $L'$ and that
$\psi({\rm supp}(f))$ strictly sits in the left side of $L'$.
Applying  Lemma~\ref{lem:Sikorav} to $K=\psi({\rm supp}(\phi))$ we
get a Hamiltonian diffeomorphism $\widetilde{g}'$ such that
$\|\widetilde{g}'\|_k= \|\widetilde{g}\|_k<E({\rm supp}(\phi))+\delta$ and that $\widetilde{g}'(K)$ lies in the right
side of $L'$ and hence
$$
\widetilde{g}'(\psi({\rm supp}(\phi)))\cap\psi({\rm supp}(f))=\emptyset.
$$
It follows that
$C_{(\psi^{-1}\widetilde{g}'\psi)\phi(\psi^{-1}\widetilde{g}'\psi)^{-1}}f=f$
since $f({\rm supp}(f))={\rm supp}(f)$. So Lemma~\ref{lem:estimate}
leads to
$$
d_k(f,C_\phi f)\leq2^{3k+2}\|\psi^{-1}\widetilde{g}'\psi\|_k=
2^{3k+2}\|\widetilde{g}'\|_k<2^{3k+2}E({\rm supp}(\phi))+2^{3k+2}\delta.
$$
Since $\delta$ can be chosen arbitrarily small,  the desired
estimate is obtained. \hfill$\Box$\vspace{2mm}

Finally, as in \cite[5.2]{LaPo97} using
Lemma~\ref{lem:biInvariance} and Proposition~\ref{th:BI(R^2n)}
we deduce that the functions
\begin{equation}\label{e:S-bi-metric}
D_k(f,g)=r_{\infty,k}(g^{-1}f),\quad k=0,1,2,\cdots
\end{equation}
give a sequence of non-degenerate bi-invariant quasi-metrics on the group
${\rm Symp}_c(\R^{2n},\omega_0)$ of all compactly supported
symplectomorphisms of $(\R^{2n},\omega_0)$.

\subsection{Banyaga's method}

Let $(M,\omega)$ be a closed (i.e., compact and without boundary) symplectic manifold.
For a smooth path $[a,b]\ni t\mapsto\phi_t\in {\rm Symp}(M,\omega)$, (which means
the mapping $(x,t)\mapsto \phi_t(x)$ to be smooth), it determines a unique smooth family
of symplectic vector fields, $\dot{\phi}_t(x)=\frac{d\phi_t}{dt}\circ\phi_t^{-1}(x)$,
whose  dual $1$-form $i_{\dot{\phi_t}}\omega$ is closed. When $[a,b]=[0,1]$ and $\phi_0={\rm id}$, $\Phi:=\{\phi_t\}$ is called a
 {\bf symplectic isotopy}  on $(M,\omega)$.

Fix a Riemannian metric $g$ on $M$.
A symplectic vector field $X$ on $M$ is said to be a {\bf harmonic vector field}
          if $i_X\omega$ is a harmonic form; and a
           smooth path $[a,b]\ni t\mapsto\phi_t\in {\rm Symp}(M,\omega)$
         is called  {\bf harmonic} if each form $i_{\dot{\phi}_t}\omega$ is harmonic.
   In particular, a
          symplectic isotopy $\Phi=\{\phi_t\}$ on $(M,\omega)$
  is called a {\bf harmonic isotopy} if it is a harmonic path in ${\rm Symp}(M,\omega)$.

  \begin{lemma}[\hbox{\cite[Lemma 1]{Ban07}}]\label{lem:IsoDecom}
Any smooth path $[a,b]\ni t\mapsto \phi_t\in {\rm Symp}(M,\omega)$  can be decomposed in a unique way as $\phi_t = \rho_t\psi_t$,  where $[a,b]\ni t\mapsto \rho_t\in {\rm Symp}(M,\omega)$ is a harmonic path and
 $[a,b]\ni t\mapsto \psi_t\in {\rm Ham}(M,\omega)$
 is a (smooth) Hamiltonian path. In particular,   if $\phi_t$ is a Hamiltonian path, then $\phi_t = \psi_t$ and $\rho_t = {\rm id}$.
\end{lemma}


  Fix a basis $\{h_1,\cdots,h_r\}$ of harmonic $1$-forms, where $r=\dim H^1(M,\mathbb{R})$.
  The space of harmonic $1$-forms on $M$ is  equipped with the following Euclidean metric:
\begin{equation}\label{def:Euclid norm}
|h|=\sum_i|\lambda_i|\quad\hbox{if}\; h=\sum_i\lambda_ih_i.
  \end{equation}
  Banyaga defined in \cite{Ban07} the length of a
symplectic isotopy $\Phi=\{\phi_t\}$  by
\begin{equation}\label{def:H-L length}
l_{HL}(\Phi)=\int_0^1\big(|\mathscr{H}_t|+(\max_x U_t-\min_x U_t)\big)dt,
\end{equation}
where $\mathscr{H}_t$ and $U_t$ are smooth families of harmonic $1$-forms and functions
respectively and satisfy the Hodge decomposition
\begin{equation}\label{e:Hodge decomp}
i_{\dot{\phi}_t}\omega=\mathscr{H}_t+dU_t.
\end{equation}
Clearly, for a Hamiltonian isotopy $\Phi$ the  formula (\ref{def:H-L length})
reduces to (\ref{def:length}).   As for (\ref{def:length}) we can prove that $l_{HL}(\Phi)$ is independent of the choice of
parametrization of the path $\Phi$. However we don't have $l_{HL}(\Phi)=l_{HL}(\Phi^{-1})$ in general, where
$\Phi^{-1}=\{\phi_t^{-1}\}$.

In \cite{Ban07}  the energy $e(\phi)$ of  any $\phi\in{\rm Symp}_0(M,\omega)$ is defined by
$$
e(\phi)=\inf_\Phi(l_{HL}(\Phi)),
$$
where $\Phi$ runs over all symplectic isotopies connecting the identity and $\phi$.
 The so called  \textbf{Hofer-like metric} is the map
\begin{eqnarray*}
\|\cdot\|_{HL}:{\rm Symp}_0(M,\omega)\rightarrow\mathbb{R}\cup\{\infty\},\;
\phi\mapsto\frac{1}{2}\big(e(\phi)+e(\phi^{-1})\big),
\end{eqnarray*}
which is actually a norm on ${\rm Symp}_0(M,\omega)$.
Notice that the norm $\|\cdot\|_{HL}$ depends on the choice of the Riemannian
metric $g$ on $M$ and the choice of the Euclidean norm $|\cdot|$ on the space
of harmonic 1-forms. However, different choices for $g$ and $|\cdot|$ yield
equivalent metrics. The {\bf Hofer-like distance} $d_{HL}$ on each connected component of
${\rm Symp}(M,\omega)$ is defined by
$d_{HL}(\phi,\psi):=\|\phi\psi^{-1}\|_{HL}$.
It is right invariant, but not left invariant.
When $\Phi$ is a Hamiltonian isotopy, (\ref{def:H-L length}) reduces to
(\ref{def:length}), thus $\|\phi\|_{HL}\leq\|\phi\|_H\;\forall\phi\in{\rm Ham}(M,\omega)$.
Moreover the subgroup ${\rm Ham}(M,\omega)$ is closed in ${\rm Symp}(M,\omega)$
endowed with the metric topology defined by $\|\cdot\|_{HL}$.
As expected by Banyaga, Buss and Leclercq \cite{BuLe} showed that
  the restriction of the Hofer-like metric to
${\rm Ham}(M,\omega)$ is equivalent to the Hofer metric.

Given any smooth symplectic path $\alpha:[a,b]\rightarrow {\rm Symp}(M,\omega)$,
we have a decomposition as (\ref{e:Hodge decomp}), $i_{\dot{\alpha}(t)}\omega=\mathscr{H}_t+dU_t$.
For every integer $k=0,1,2,\cdots$, we define the {\bf $k$-length} of $\alpha$ as
\begin{equation}\label{def:H-L k}
l_{HL,k}(\alpha):=\sum_{i=0}^{k}\int_a^b\left[\left|\frac{\partial^i\mathscr{H}_t}
{\partial t^i}\right|+\left(\max_x\frac{\partial^iU_t}{\partial t^i}-
\min_x\frac{\partial^iU_t}{\partial t^i}\right)\right]dt.
\end{equation}
where $|\cdot|$ is as in (\ref{def:Euclid norm}).
Obviously, $l_{HL,0}(\alpha)=l_{HL}(\alpha)$. But when $k\geq1$, $l_{HL,k}$
depends on the choice of parametrization of the path which is different from
$l_{HL}$.

 A continuous path
$\Phi:[0,1]\rightarrow{\rm Symp}_0(M,\omega)$ is called a {\bf piecewise smooth
symplectic isotopy} if there exists a division
$0=t_0<t_1<\cdots<t_n=1, n\in\mathbb{N}$, such that for each $i=1,\cdots,n$,
$\Phi_i=\Phi|_{[t_{i-1},t_i]}$ is smooth, and $\Phi(0)={\rm id}$.
We define the {\bf $k$-length} of $\Phi$ as
$$
l_{HL,k}(\Phi):=\sum_{i=1}^nl_{HL,k}(\Phi_i).
$$
For any $\phi\in {\rm Symp}_0(M,\omega)$, let
$\Omega(\phi)$ be the set of all piecewise smooth isotopies
$\Phi:[0, 1]\to {\rm
Symp}_0(M,\omega)$ with  $\Phi(1)=\phi$.
Define the energy of $\phi$ by
$$
e_k(\phi):=\inf\left\{l_{HL,k}(\Phi)\,\Big|\,\Phi\in\Omega(\phi)\right\},
$$
and the corresponding {\bf Hofer-like $k$-metric} by
$$
\|\phi\|_{HL,k}=\frac{1}{2}\big(e_k(\phi)+e_k(\phi^{-1})\big),
\quad k=0,1,\cdots.
$$
We define the {\bf Hofer-like $k$-distance}  by
$$
d_{HL,k}(\phi,\psi):= \|\phi\psi^{-1}\|_{HL,k}
$$
on every connected component of ${\rm Symp}(M,\omega)$.

\begin{prop}\label{prp:property2}
On a closed symplectic manifold $(M, \omega)$ with a fixed Riemannian metric $g$, we have
\begin{description}
\item[(i)] $\|\cdot\|_{HL}=\|\cdot\|_{HL,0}\leq\|\cdot\|_{HL,1}
\leq\|\cdot\|_{HL,2}\leq\cdots$,
\item[(ii)] for each $k$, $\|\cdot\|_{HL,k}$ is a quasi-norm on
${\rm Symp}_0(M,\omega)$.
\end{description}
\end{prop}

\noindent{\bf Proof.} (i) Notice that the corresponding Lemma~\ref{lem:HoferM} still hold
in the current context. Hence for $\phi\in{\rm Symp}_0(M,\omega)$, we have
$e_0(\phi)=e(\phi)$, so $\|\cdot\|_{HL,0}=\|\cdot\|_{HL}$.

\noindent{(ii)} \;\;\underline{Symmetry}: By definition, if $\phi,\psi\in
{\rm Symp}(M,\omega)$ is in the same connect component, then
$d_{HL,k}(\phi,\psi)=d_{HL,k}(\psi,\phi)$.

\underline{Non-degeneracy}: By (i) and the non-degeneracy of $\|\cdot\|_{HL}$,
this is obvious.

\underline{The quasi-triangle inequality}: Since the proof is similar to the one
of property (ii) in Theorem~\ref{th:property}, we skip some details and only
outline the ideas of the proof. For any two symplectomorphisms
$\phi,\psi\in {\rm Symp}_0(M,\omega)$, we choose
\begin{eqnarray*}
&&\Omega(\phi)\ni\Phi:[0,1] \rightarrow{\rm Symp}_0(M,\omega)\quad\hbox{and}\\
&&\Omega(\psi)\ni\Psi:[0,1] \rightarrow{\rm Symp}_0(M,\omega).
\end{eqnarray*}
Define the concatenation product $\Phi\ast\Psi$ of $\Phi$ and $\Psi$ by
$$
\Phi\ast\Psi(t)=
\begin{cases}
\Psi(2t) & 0\leq t\leq 1/2\\
\Phi(2t-1)\psi & 1/2\leq t\leq 1.
\end{cases}
$$
then $\Phi\ast\Psi\in\Omega(\phi\psi)$. By definition we have
\begin{eqnarray*}
e_k(\phi\psi)&\leq &l_{HL,k}(\Phi\ast\Psi)
\leq2^k(l_{HL,k}(\Psi)+
l_{HL,k}(\Phi(t)\psi))\\
&=&2^k(l_{HL,k}(\Phi)+l_{HL,k}(\Psi))
\end{eqnarray*}
 for all $\Phi\in\Omega(\phi)$, $\Psi\in\Omega(\psi)$.
Taking the infimum respectively we get
\begin{eqnarray*}
e_k(\phi\psi)&\leq &2^k(\inf\{l_{HL,k}(\Phi)|\Phi\in\Omega(\phi)\}+
\inf\{l_{HL,k}(\Psi)|\Psi\in\Omega(\psi)\})\\&=&2^k(e_k(\phi)+e_k(\psi)).
\end{eqnarray*}
That is,  $e_k$ (and so $\|\cdot\|_{HL,k}$) satisfies the quasi-triangle inequality.
\hfill$\Box$\vspace{2mm}

\begin{prop}\label{prp:property3}
Let $(M,\omega)$ be a closed symplectic manifold. Then for each $k$
the subgroup ${\rm Ham}(M,\omega)$ is closed in ${\rm Symp}(M,\omega)$
with respect to the metric topology defined by $\|\cdot\|_{HL,k}$.
\end{prop}

\noindent{\bf Proof.}
The ideas are  similar to those of Theorem~14.2.A in \cite{PoBook}.
 Suppose there exists
a sequence $\{f_n\}\subset{\rm Ham}(M,\omega)$ and $\phi\in{\rm Symp}(M,\omega)$,
satisfying $d_{HL,k}(f_n,\phi)\rightarrow0$ when $n\rightarrow\infty$. We intend
to prove $\phi\in{\rm Ham}(M,\omega)$.

Since $\lim_{n\rightarrow\infty}\|f_n\phi^{-1}\|_{HL,k}=0$, by the definition of
$\|\cdot\|_{HL,k}$, for $\forall \varepsilon>0$, $\exists N_0>0$, such that for each $N\geq N_0$
we have
$$
\inf\{l_{HL,k}(\Phi)|\Phi\in\Omega(f_N\phi^{-1})\}<\varepsilon,
$$
and so a $\Phi^N\in\Omega(f_N\phi^{-1})$  such that
$l_{HL}(\Phi^N)\le l_{HL,k}(\Phi^N)<\varepsilon$.

Assume the division of $\Phi^N$ is given by $0=t_0<t_1<\cdots<t_n=1$. For each
$i=1,\cdots,n$, $\Phi^N_i=\Phi^N|_{[t_{i-1},t_i]}$ is a smooth symplectic path.
Next we translate $\Phi^N$ into a smooth symplectic isotopy through the precedure
of reparametrization. We choose a increasing, surjective smooth function
$s_i:[t_{i-1},t_i]\rightarrow[t_{i-1},t_i]$ for each $i$, and require $s_i$
is constant near the ends of its interval of definition. Define
$\widetilde{\Phi}^N:[0,1]\rightarrow {\rm Symp}_0(M,\omega)$ by
$$
\widetilde{\Phi}^N(t) =\Phi^N_i(s_i(t)),\quad\forall
t\in[t_{i-1},t_i],\;i=1,\cdots,n,
$$
then $\widetilde{\Phi}^N$ is a smooth symplectic isotopy. In fact if the
harmonic $1$-forms and Hamiltonian functions generated by $\Phi^N$ are
$\{^i\mathscr{H}^N_t\}_{i=1}^n$ and $\{^iU^N_t\}_{i=1}^n$ respectively,
and the harmonic $1$-forms and Hamiltonian functions generated by $\widetilde{\Phi}^N$
are $\widetilde{\mathscr{H}}^N_t$ and $\widetilde{U}^N_t$ respectively, then
when $t\in[t_{i-1},t_i]$, we have
\begin{equation}\label{e:reparametrization}
\widetilde{\mathscr{H}}^N_t=s_i'(t)\cdot^i\mathscr{H}^N_{s_i(t)}\;\;,\qquad
\widetilde{U}^N_t=s_i'(t)\cdot^iU^N_{s_i(t)}.
\end{equation}
By the change of variable formula, we get
$l_{HL}(\widetilde{\Phi}^N)=l_{HL}(\Phi^N)$,
In particular,
$$
\int_0^1\left|\widetilde{\mathscr{H}}^N_t\right|dt=
\sum_{i=1}^n\int_{t_{i-1}}^{t_i}\left|^i\mathscr{H}^N_t\right|dt.
$$
Since $l_{HL}(\Phi^N)<\varepsilon$, we get
$\int_0^1|\widetilde{\mathscr{H}}^N_t|dt<\varepsilon$.

Recall that the flux is a surjective homomorphism from
the universal covering space $\widetilde{{\rm Symp}_0}(M,\omega)$  of
${\rm Symp}_0(M,\omega)$ to $H^1(M,\mathbb{R})$ given by
\begin{equation}\label{def:flux}
{\rm Flux}(\{\phi_t\})=\left[\int_0^1(i_{\dot{\phi}_t}\omega)dt\right]
\in H^1(M,\mathbb{R}).
\end{equation}
$\Gamma_\omega:={\rm Flux}\bigl(\pi_1({\rm Symp}_0(M,\omega))\bigr)\subset H^1(M,\mathbb{R})$
is called the {\bf flux group} of $(M,\omega)$, and is discrete as proved by Ono  in \cite{Ono}.
Flux descends to a surjective homomorphism
$$
{\rm flux}:{\rm Symp}_0(M,\omega)\rightarrow H^1(M,\mathbb{R})/\Gamma_\omega
$$
with kernel ${\rm Ham}(M,\omega)$ (cf. \cite{McSa}).

For any symplectic isotopy $\Phi$, let $\mathscr{H}(\Phi)$ denote the harmonic
representation of the cohomology class ${\rm Flux}(\Phi)$. The decomposition
(\ref{e:Hodge decomp}) implies that
$$
\mathscr{H}(\Phi)= \int_0^1\mathscr{H}_tdt.
$$
It follows from this that
\begin{eqnarray*}\label{e:harmonic estimate}
\left|\mathscr{H}(\widetilde{\Phi}^N)\right|
&=&\left|\int_0^1\widetilde{\mathscr{H}}^N_tdt \right|
=\left|\int_0^1\sum_{i=1}^r\widetilde{\lambda}^N_i(t)h_idt\right|\nonumber\\
&=&\sum_{i=1}^r\left|\int_0^1\widetilde{\lambda}^N_i(t)dt\right|\nonumber\\
&\leq&\int_0^1\sum_{i=1}^r|\widetilde{\lambda}^N_i(t)|dt
=\int_0^1|\widetilde{\mathscr{H}}^N_t|dt<\varepsilon,
\end{eqnarray*}
where $\widetilde{\mathscr{H}}^N_t$ is decomposed as
$\sum_{i=1}^r\widetilde{\lambda}^N_i(t)h_i$.

Starting from $\widetilde{\Phi}^N$, we could construct a smooth symplectic
path $\widetilde{\Phi}^N\circ\phi$ connecting $\phi$ and $f_N$. Obviously we have
$\mathscr{H}(\widetilde{\Phi}^N)=\mathscr{H}(\widetilde{\Phi}^N\circ\phi)$.
Choose any symplectic isotopy $\Psi$ from ${\rm id}$ to $\phi$, and any
Hamiltonian isotopy $\alpha^N$ from ${\rm id}$ and $f_N$, we get a loop  $(-\alpha^N)\sharp(\widetilde{\Phi}^N\circ\phi)\sharp\Psi$, whose flux has
 the harmonic representation
$$\mathscr{H}(\Psi)+
\mathscr{H}(\widetilde{\Phi}^N\circ\phi)\in\Gamma_\omega
$$
because  a Hamiltonian path has zero flux.
Note that $|\mathscr{H} (\widetilde{\Phi}^N)|<\varepsilon$, and that $\varepsilon$
is arbitrary small. We deduce that $\mathscr{H}(\Psi)\in\Gamma_\omega$
since $\Gamma_\omega$ is discrete. Hence $\phi\in\ker({\rm flux})={\rm Ham}(M,\omega)$. \hfill$\Box$\vspace{2mm}

\begin{rmk}
{\rm By the method in \cite{Ban07}, we can't obtain that $f_N^{-1}\phi$ is a
Hamiltonian diffeomorphism for every $N$ large enough, but could only get
the distance from ${\rm flux}(f_N^{-1}\phi)$ to $\Gamma_\omega$ trends to zero as $N\rightarrow\infty$. }
\end{rmk}

Han \cite{Han} also  introduced a method constructing bi-invariant
(quasi) metrics on ${\rm Symp}(M,\omega)$ from the Hofer metric.
For a fixed positive number $K$, he  defined
$\|\phi\|_K=\min(\|\phi\|_H, K)$ if $\phi\in{\rm Ham}(M, \omega)$,
and $K$ otherwise. However,  when the above defined quasi-metrics (or metrics) $r_{\alpha,k}$, $\|\cdot\|_K$
are restricted back to ${\rm Ham}(M, \omega)$, the induced topologies are
in general different from that of the Hofer metric.

\section{Concluding remarks} \setcounter{equation}{0}

 Extensions of the Hofer metric to contact geometry were also studied, see
  Banyaga and  Donato \cite{Ban06},  Banyaga and Spaeth \cite{BS08} and
   M\"{u}ller and Spaeth \cite{MuSp}.
Our proceeding constructions can be completed in contact manifolds.
Let $(N,\alpha)$ be a compact contact manifold of dimension $2n+1$.
There exists a one-to-one correspondence between contact isotopies on $(N,\alpha)$
and elements of the space $C^\infty(N\times[0,1])$, $\{f_t\}\leftrightarrow H$,
where $i_{X_t}\alpha=H_t$ with $H_t=H(\cdot,t)$ and $X_t=(\frac{d}{dt}f_t)\circ f_t^{-1}$;
$H$ is called the contact Hamiltonian function of $\{f_t\}$.
Call $\phi\in{\rm Diff}(N)$ a strictly contact diffeomorphism if $\phi^\ast\alpha=\alpha$.
A contact isotopy is said to be strictly if each contact diffeomorphism in the isotopy
is strictly. Denote by $G_\alpha(N)$ the group of strictly contact diffeomorphisms
which are strictly contact isotopic to the identity.

Consider
a surjective homomorphism from the universal cover $\widetilde{G_\alpha(N)}$ of
$G_\alpha(N)$ to $\mathbb{R}$ given by
$$
\{\phi_t\}\mapsto c(\{\phi_t\})=\frac{1}{{\rm Vol}(N)}\int_0^1\Big(\int_N
H_t(x)\nu_\alpha\Big)dt,
$$
where $H_t$ is the contact Hamiltonian of $\{\phi_t\}$ and
the canonical volume form
$\nu_\alpha:=\alpha\wedge(d\alpha)^n$. For each $k=0,1,2,\cdots$, we define
{\bf $k$-contact length} of $\phi_t$ by
\begin{eqnarray*}
{\rm Length}_{c,k}(\{\phi_t\}):=|c(\phi_t)|+\sum_{i=0}^k\int_0^1
\left(\max_{x\in N}\frac{\partial^i H}{\partial t^i}(x,t)-\min_{x\in
N} \frac{\partial^i H}{\partial t^i}(x,t)\right)dt,
\end{eqnarray*}
and {\bf $k$-contact energy} of $\phi\in G_\alpha(N)$ by
$$
E_{c,k}(\phi)=\inf_{\phi_t}\big({\rm Length}_{c,k}(\{\phi_t\})\big),
$$
where $\{\phi_t\}$ takes over all piecewise smooth strictly contact isotopy from ${\rm id}$ to $\phi$.
If $k=0$ it becomes the contact length and contact energy in \cite[(17) and (20)]{Ban06}.
Using the results in  \cite{Ban06, MuSp} we may directly prove

\begin{thm}
For each $k=0,1,\cdots$, the mapping
\begin{eqnarray*}
d_{c,k}:G_\alpha(N)\times G_\alpha(N)\rightarrow [0,\infty),\quad
(\phi,\psi)\mapsto E_{c,k}(\phi\psi^{-1})
\end{eqnarray*}
 is a bi-invariant quasimetric on $G_\alpha(N)$.
\end{thm}

When the contact manifold $(N,\alpha)$ is regular, that is, the Reeb field $R_\alpha$ of $\alpha$
generates a free $S^1$-action on $N$, the quotient manifold $B=N/S^1$ is a base of a principal $S^1$-
bundle $\pi:N\to B$ and $B$ has a canonical symplectic form $\omega$ satisfying $\pi^\ast\omega=d\alpha$.
In this case there exists an exact sequence
$$\{1\}\to S^1\to G_\alpha(M)\stackrel{p}{\longrightarrow} {\rm Ham}(B,\omega)\to\{1\}.
$$
As in the proof of \cite[Lemma 4.2]{BS08} it is not hard to prove that
$E_{c,k}(\phi)\ge\|p(\phi)\|_k$ for
any $\phi\in G_\alpha(N)$.

As in Hofer geometry it is an important topic to study  geodesics of our metrics.

\appendix
\section{Appendix: Semigroupoid Metrization Theorem} \setcounter{equation}{0}

Given a semigroupoid $(G,\ast)$,
let $G^{(1)}=G$, $G^{(2)}=\{(a,b)\in G\times G: a\ast b\;\hbox{is well-defined}\}$
and for each $N\in\N$, $N\ge 2$ let
\begin{eqnarray*}
G^{(N)}:=\bigl\{(a_1,\cdots,a_N)\in G\times\cdots\times G &|& (a_j, a_{j+1})\in G^{(2)}\\
&& \forall j\in\{1,\cdots,N-1\}\bigr\}.
\end{eqnarray*}
In particular, if  $(G,\ast)$ is a semigroup, $G^{(N)}$ is just the Cartesian product of $N$ copies of $G$.

\begin{theorem}{\rm (\hbox{\cite[Cor.3.33]{MiMMM}})}\label{th:MiMMM}
Let $(G,\ast)$ be a semigroupoid, and assume that $\psi:G\to [0,\infty]$ is a function
with the property that there exists a finite
constant $C\geq1$  such that
\begin{equation}\label{e:MiMMM1}
\psi(a\ast b)\le C(\psi(a)+\psi(b)),\quad\hbox{for all}\;(a,b)\in
G^{(2)}.
\end{equation}
Introduce
\begin{equation}\label{e:MiMMM2}
\alpha:=\frac{1}{1+\log_2C}\in (0, 1]
\end{equation}
and define the function $\psi_\sharp:G\to [0, \infty]$ by
\begin{equation}\label{e:MiMMM3}
\psi_\sharp(a):=\inf\left\{\left(\sum^N_{i=1}\psi(a_i)^\alpha\right)^{\frac{1}{\alpha}}:\;N\in\N,\;
(a_1,\cdots,a_N)\in G^{(N)},\;a=a_1\ast\cdots\ast a_N\right\}.
\end{equation}
Then $\psi\approx\psi_\sharp$. More specifically, with $C$ the same
constant as in (\ref{e:MiMMM1}), one has
\begin{equation}\label{e:MiMMM4}
(2C)^{-2}\psi\le\psi_\sharp\le\psi\quad\hbox{on}\;G.
\end{equation}
In particular, $\psi^{-1}(\{0\})=\psi_\sharp^{-1}(\{0\})$.
 Furthermore, for every $\beta\in (0, \alpha]$ one has
\begin{equation}\label{e:MiMMM5}
\psi_\sharp(a\ast
b)^\beta\le\psi_\sharp(a)^\beta+\psi_\sharp(b)^\beta,\quad\forall
(a,b)\in G^{(2)},
\end{equation}
 and $\psi_\sharp=\psi$ on $G$ if and
only if $\psi_\sharp(a\ast
b)^\alpha\le\psi_\sharp(a)^\alpha+\psi_\sharp(b)^\alpha$ for all
$(a,b)\in G^{(2)}$. Finally, for each $N\in\N$ the original function
$\psi$ satisfies
\begin{equation}\label{e:MiMMM6}
\psi(a_1\ast\cdots a_N) \le
4C^2\left\{\sum^N_{i=1}\psi(a_i)^\beta\right\}^{\frac{1}{\beta}}
\end{equation}
whenever $a_1,\cdots,a_N\in G$ are such that
\begin{equation}\label{e:MiMMM7}
(a_i, a_{i+1})\in G^{(2)}\quad\hbox{for every}\;i\in\{1,\cdots,
N-1\}.
\end{equation}
In particular, if $(a_i)_{i\in\N}\subset G$  is a sequence with the
property that (\ref{e:MiMMM7}) holds for every number $N\in\N$ with
$N\ge 2$, then for each finite number $\beta\in (0, \alpha]$ one has
\begin{equation}\label{e:MiMMM8}
\sup_{N\in\N}\psi(a_1\ast\cdots a_N) \le
4C^2\left\{\sum^\infty_{i=1}\psi(a_i)^\beta\right\}^{\frac{1}{\beta}}.
\end{equation}
\end{theorem}

\begin{remark}\label{rm:MiMMM}
{\rm When $(G,\ast)$ is a group, a function $\psi:G\to
[0,\infty]$ is said to be {\bf conjugate invariant} provided
 $\psi(b\ast a\ast b^{-1})=\psi(a)\;\forall
a,b\in G$.  In this case $\psi_\sharp$ is also conjugate invariant. In fact, for any $b\in G$,
$(a_1,\cdots,a_N)\in G^{(N)},\;a=a_1\ast\cdots\ast a_N$, since
$b\ast a\ast b^{-1}=(b\ast a_1\ast b^{-1})\cdots\ast (b\ast a_N\ast
b^{-1})$ and
$$
\psi(b\ast a_i\ast b^{-1})=\psi(a_i)\quad\forall i=1,\cdots,N,
$$
it follows from the definition of $\psi_\sharp(a)$ in
$(\ref{e:MiMMM3})$ that $\psi_\sharp(b\ast a\ast
b^{-1})=\psi_\sharp(a)$.}
\end{remark}


\end{document}